\documentclass[11pt]{amsart}

\usepackage{mathrsfs}
\usepackage[mathcal]{euscript}
\usepackage{graphicx}
\usepackage{amsthm,amscd}
\usepackage{calc}
\usepackage{amsmath,amssymb,amsfonts}
\usepackage{enumerate}
\usepackage{indentfirst}
\usepackage[all]{xy}
\usepackage[
  a4paper,
  includeheadfoot,
  margin=25mm,
  headheight=14pt,
  headsep=7mm,
  footskip=10mm
]{geometry}
\usepackage[pagebackref]{hyperref}

\renewcommand*{\backref}[1]{}

\renewcommand*{\backrefalt}[4]{%
  \ifcase #1\relax
  \or\space#2%
  \else\space#2%
  \fi}

\numberwithin{equation}{section}

\theoremstyle{plain}
\newtheorem{thm}{Theorem}[section]
\newtheorem{lemma}[thm]{Lemma}

\newtheorem{prop}[thm]{Proposition}
\newtheorem{cor}[thm]{Corollary}

\theoremstyle{definition}

\newtheorem{defn}[thm]{Definition}
\newtheorem{rem}[thm]{Remark}
\newtheorem{notation}[thm]{Notation}
\newtheorem{step}{Step}

\newcommand{\C}{\mathbb C}
\newcommand{\Z}{\mathbb Z}
\newcommand{\Q}{\mathbb Q}
\newcommand{\R}{\mathbb R}
\newcommand{\cC}{\mathcal C}

\newcommand{\cE}{\mathcal E}
\newcommand{\cF}{\mathcal F}

\newcommand{\cL}{\mathcal L}

\newcommand{\cN}{\mathcal N}
\newcommand{\cO}{\mathcal O}

\newcommand{\cS}{\mathcal S}
\newcommand{\cX}{\mathcal X}
\newcommand{\Pic}{\operatorname{Pic}}

\newcommand{\Supp}{\operatorname{Supp}}
\newcommand{\Vol}{\operatorname{Vol}}

\begin{document}

\title[Upper semicontinuity of algebraic dimension]
{Upper Semicontinuity of Algebraic Dimension}

\author{Sheng Rao}
\address{School of Mathematics and Statistics, Wuhan University, Wuhan 430072, China}
\email{likeanyone@whu.edu.cn}
\thanks{The author is  partially supported by NSFC (Grant No. 12271412, W2441003, 12671107) and Hubei Provincial Innovation Research Group Project (Grant No. 2025AFA044).}

\date{}



\begin{abstract}
We prove the conjecture that the deformation limit of Moishezon manifolds under a smooth deformation over the unit disc in $\mathbb{C}$ is Moishezon.
More generally, for a proper holomorphic submersion over the unit disc
with connected compact fibers satisfying the $\partial\bar\partial$-lemma
except possibly at one parameter, we show that the very general algebraic
dimension is the minimum on the entire disc.
The proof builds on Barlet's theory of cycle spaces and 
Rao--Tsai's torsion-freeness method for extending integral Chern classes. Taking this approach as a starting point, we handle the
remaining finite-order torsion obstruction by a scalar phase estimate
for divisor masses.  This gives uniform Gauduchon mass bounds on each
component over compact subsets of the base without assuming
torsion-freeness.  Bishop compactness then establishes properness of
every irreducible component of the relative divisor space, and Barlet's
algebraic-dimension theorem yields the conclusion directly.
Combined with Rao--Tsai's theorem, the deformation-limit result also
shows that every fiber is Moishezon whenever uncountably many fibers
are Moishezon.
The Fujiki--Pontecorvo families and Huybrechts' twistor and Brauer families
clarify the scope and necessity of the hypotheses.
\end{abstract}

\maketitle

\setcounter{tocdepth}{1}
\tableofcontents

\section{Introduction}\label{s:introduction}

Let $\Delta:=\{t\in\C:|t|<1\}$ be the \emph{unit disc} and write
$\Delta^*:=\Delta\setminus\{0\}$.  Let
\begin{equation}\label{eq:family}
 \pi:\cX\longrightarrow\Delta
\end{equation}
be a proper holomorphic submersion with connected compact fibers
$X_t:=\pi^{-1}(t)$ of complex dimension $n$.  The \emph{algebraic dimension}
$a(X)$ of a compact connected complex manifold is the transcendence degree
over $\C$ of its field of meromorphic functions.  It satisfies
$0\leq a(X)\leq\dim_\C X$; equality on the right characterizes \emph{Moishezon
manifolds}.

According to Kodaira's classification of surfaces (e.g., \cite{Si83,Bu99,Lm99}), every compact complex surface with even first Betti number is K\"ahler. Thus Ehresmann's theorem \cite{Eh47} and the deformation invariance of Hodge numbers imply that the limit surface $X_{0}$ of compact complex surfaces in Fujiki class $\mathcal{C}$ (or equivalently K\"ahler surfaces) is K\"ahler, while \cite[Corollary 1.3]{Po13} proves that the deformation limit of Moishezon (or equivalently projective) surfaces is projective.

In the higher-dimensional setting, Hironaka's counterexample \cite{Hi62} demonstrates that even if all fibers \( X_t \) are projective for \( t \in \Delta^* \), the central fiber \( X_0 \) need not be projective or K\"ahler.

D. Popovici \cite{Po13} proved that if all fibers \( X_t \) are projective for \( t \in \Delta^* \) and the central fiber \( X_0 \) either satisfies the deformation invariance of the Hodge number \( h^{0,1} \) or admits a strongly Gauduchon metric, then \( X_0 \) is Moishezon.
A Hermitian metric $\omega$ on a compact complex $n$-dimensional manifold
is called a \emph{Gauduchon metric} if
$\partial\bar\partial\omega^{n-1}=0$, and a
\emph{strongly Gauduchon metric} if $\partial\omega^{n-1}$ is
$\bar\partial$-exact, as in \cite{Po13}.
The existence of Gauduchon metrics on any compact complex manifold was established in \cite{Ga77}. Moreover, if \( X \) satisfies the \( \partial \bar{\partial} \)-lemma, then every Gauduchon metric on \( X \) is automatically strongly Gauduchon. Here the \emph{$\partial\bar\partial$-lemma} means
the following: for every pure-type $d$-closed form on a complex manifold, the properties of
$d$-exactness, $\partial$-exactness, $\bar\partial$-exactness and
$\partial\bar\partial$-exactness are equivalent. 

D. Barlet \cite{Ba15} further extended these results, showing that if all fibers \( X_t \) are Moishezon for \( t \in \Delta^* \) and \( X_0 \) satisfies either of the aforementioned conditions, then \( X_0 \) remains Moishezon.  
More precisely, his two semicontinuity results can be stated together as
follows.

\begin{thm}[{\cite[Theorems~1.0.2 and~1.0.3]{Ba15}}]
\label{thm:barlet-semicontinuity}
Let $\pi:\cX\to\Delta$ be the family in \eqref{eq:family}.
\begin{enumerate}[$(i)$]
\item If $X_{t_0}$ admits a smooth strongly Gauduchon form, there are an
open neighbourhood $U\subset\Delta$ of $t_0$, an integer $a\ge0$, and a
countable union
$\Sigma$ of irreducible proper analytic subsets of $U$ 
such that
\[
 a(X_t)\ge a\quad(t\in U),
 \qquad a(X_t)=a\quad(t\in U\setminus\Sigma).
\]
\item Suppose that $\Delta^\circ\subset\Delta$ is a dense analytic Zariski
open subset and that every $X_t$ with $t\in\Delta^\circ$ satisfies the
$\partial\bar\partial$-lemma.
If $h^{0,1}(X_t)$ is constant on $\Delta$, then
\[
 a(X_t)\ge\inf_{u\in\Delta^\circ}a(X_u)\qquad(t\in\Delta).
\]
\end{enumerate}
\end{thm}

More recently, S. Rao--I.-H. Tsai obtained the following result.

\begin{thm}[{\cite[Theorem~1.4]{RT21}}]
\label{thm:rao-tsai}
Let $\pi:\cX\to\Delta$ be the family in \eqref{eq:family}, and assume
that uncountably many fibers are Moishezon.
For any $t_0\in\Delta$, the fiber $X_{t_0}$ is Moishezon if either
\begin{enumerate}[$(i)$]
\item $R^2\pi_*\cO_{\cX}$ is torsion-free in a neighbourhood of $t_0$
(equivalently, $h^{0,1}(X_t)$ is constant in a neighbourhood of $t_0$,
by \cite[Proposition~4.12]{RT21}); or
\item $X_{t_0}$ admits a strongly Gauduchon metric.
\end{enumerate}
\end{thm}

J. Koll\'ar \cite{Ko22} proposed a conjecture in the case of flat families, asserting that if all fibers \( X_t \) are Moishezon for \( t \in \Delta^* \) and the central fiber \( X_0 \) is irreducible with rational singularities in a proper flat morphism over \( \Delta \), then \( X_0 \) is Moishezon. A brief review of progress on this topic, along with a collection of further stimulating open questions, can be found in \cite[Section 1]{Ko22}. 

The main result of this paper is the following.

\begin{thm}[Upper semicontinuity of algebraic dimension]\label{thm:main}
Assume that there is a point $t_0\in\Delta$ such that $X_t$ satisfies
the $\partial\bar\partial$-lemma for every $t\ne t_0$.
Then there are an integer $a_{\mathrm{gen}}\in\{0,\ldots,n\}$ and a
countable union $\Sigma$ of proper closed analytic subsets of $\Delta$
such that
\[
 a(X_t)=a_{\mathrm{gen}}\quad(t\in\Delta\setminus\Sigma),
 \qquad a(X_t)\ge a_{\mathrm{gen}}\quad(t\in\Delta).
\]
Thus the very general algebraic dimension is the minimum on the entire
disc.  In particular, for every $s\in\Delta$ and every sufficiently small
$\varepsilon>0$ with $\{t:|t-s|<\varepsilon\}\Subset\Delta$,
\begin{equation}\label{eq:main-inequality}
 a(X_s)\geq
 \inf_{0<|t-s|<\varepsilon}a(X_t)=a_{\mathrm{gen}}.
\end{equation}
\end{thm}

The exceptional set $\Sigma$ is countable, but need not be closed in
the analytic Zariski topology.  Equivalently, each superlevel set
$\{t\in\Delta:a(X_t)\ge k\}$ is a countable union of closed analytic
subsets; it need not itself be a closed analytic subset.  This is not
upper semicontinuity for the ordinary analytic Zariski topology.

Every Moishezon manifold satisfies the $\partial\bar\partial$-lemma
(see \cite{Pa66} and \cite[Corollary~5.23]{DGMS75}), and
$a(X_t)=n$ for $t\ne0$.  Theorem~\ref{thm:main} gives
$a(X_0)\ge n$, hence equality.
\begin{cor}[Moishezon deformation limit]\label{cor:moishezon}
If $X_t$ is Moishezon for every $t\in\Delta^*$, then $X_0$ is
Moishezon.
\end{cor}

Actually, we can obtain a more general result: 
\begin{thm}[Uncountably many Moishezon fibers]\label{Moishezon}
If uncountably many fibers of \eqref{eq:family} are Moishezon, then every
fiber is Moishezon.
\end{thm}

We outline the proof of Theorem~\ref{thm:main} in
Section~\ref{s:componentwise-main-proof}.  After an automorphism of the
unit disc, assume that the possibly exceptional parameter is $0$.
The aim is to prove properness of every irreducible component of the
relative divisor space.  Vertical components are compact, so we fix a
horizontal component $\Lambda$ and its integral divisor class
$c\in H^2(M,\Z)$, where $M$ is the common differentiable fiber supplied
by an Ehresmann trivialization.  The $\partial\bar\partial$-lemma and
base change imply that its holomorphic obstruction
$o(c)\in H^0(\Delta,R^2\pi_*\cO_{\cX})$ vanishes on $\Delta^*$.
Local total-space line bundles of class $c$ then control divisor masses,
and Bishop compactness and Remmert's theorem show that
$\Lambda^*\to\Delta^*$ is proper and surjective.  This use of fixed
integral classes builds on Rao--Tsai's torsion-freeness method
\cite[Propositions~4.12 and~4.13]{RT21}: when the direct-image sheaf is
torsion-free, local vanishing of $o(c)$ extends globally and yields a
total-space line bundle.  Without this hypothesis, the remaining
obstruction is supported at $0$ and is annihilated locally by $t^N$
for some integer $N\ge0$, by Lemma~\ref{lem:revised-discrete-support}.

The additional step is a scalar phase estimate that handles this
finite-order torsion.  For a smooth closed representative $\alpha$ of
$c$, the relation $t^No(c)=0$ yields a smooth absolute $(0,1)$-form $u$
with $\bar\partial_{\cX}u=t^N\alpha^{0,2}_{\cX}$ near the central fiber,
by Lemma~\ref{lem:componentwise-torsion-primitive}.
For a smooth family of Gauduchon forms $G_t$, this gives the mass formula
\begin{equation*}
 \Vol_{G_t}(Z)=B(t)+2\operatorname{Re}\!\left(t^{-N}A(t)\right)\ge0,
 \qquad Z\in\Lambda_t,\quad t\ne0,
\end{equation*}
where $A$ is complex-valued and $B$ is real-valued, both smooth across
$0$.  Surjectivity over $\Delta^*$ supplies a divisor for every nonzero
parameter, so nonnegativity of the masses gives a lower bound for
$\operatorname{Re}(t^{-N}A(t))$ in every angular direction.
The phase Lemma~\ref{lem:componentwise-phase}, proved by finite Taylor expansion and Fourier
orthogonality, forces $A(t)=O(|t|^N)$.  Thus divisor masses stay bounded
near $0$; density of $\Lambda^*$ and continuity of cycle integration
extend the bound to the central cycles, as established in
Proposition~\ref{prop:componentwise-mass-bound}.  Bishop compactness, in the form of
Lemma~\ref{lem:revised-bishop-cycle}, now gives
properness of $\Lambda$ over the whole disc.  Applying Barlet's
algebraic-dimension Theorem~\ref{thm:revised-barlet-local} to this
componentwise properness proves Theorem~\ref{thm:main}, without a
good-filling assumption or a separate comparison of relative Kodaira
images.  The estimate controls only the scalar Gauduchon pairing needed
for compactness, not a complete minimal $\bar\partial$-solution.

The examples in Section~\ref{s:scope-examples} distinguish two notions of
semicontinuity.  A. Fujiki--M. Pontecorvo construct a non-K\"ahler family for
which both the ordinary limsup inequality and the punctured-neighbourhood
infimum inequality fail.  Twistor and Brauer families of K3 surfaces show
that the ordinary limsup inequality may fail even when every fiber is
K\"ahler, while the weaker infimum inequality remains valid.  The two
line bundles isolated by Fujiki--Pontecorvo also show why extension of a
prescribed punctured representative is stronger than extension of its
Chern class.

The paper is organized in four sections and an appendix.
Section~\ref{s:revised-preliminaries} organizes the tools for this proof:
cycle-space compactness, Gauduchon mass pairings, and the integral-class
obstruction.  The complementary good-filling constructions are collected
separately in Section~\ref{ss:revised-good-fillings}.
Section~\ref{s:componentwise-main-proof} establishes componentwise
properness by phase estimates and proves Theorems~\ref{thm:main}
and~\ref{Moishezon}.  Section~\ref{s:scope-examples} compares the conclusion with known test
families.  Appendix~\ref{s:applications} develops applications to
Moishezon families.  Standard foundational theorems are stated with the hypotheses
used here, and the steps connecting them to the argument are explained.

\begin{notation}\label{not:conventions}
Complex spaces are reduced unless otherwise specified; Douady spaces
and the subspaces they classify may be nonreduced.  Cycles retain their
positive integral multiplicities.  The base disc may be shrunk without
further mention.
After choosing an Ehresmann trivialization, all fiber cohomology groups are
identified with those of a fixed differentiable manifold $M$.  A subscript
$t$ denotes restriction to $X_t$.  Integration currents are normalized so
that the current $[D]$ of an effective Cartier divisor represents the real
image of $c_1(\cO(D))$.
For a coherent analytic sheaf $\cF$ on a compact complex space $Y$, we write
\[
 h^q(Y,\cF):=\dim_{\C}H^q(Y,\cF).
\]
In particular, $h^q(Y_b,\cF|_{Y_b})$ denotes the complex dimension of
the cohomology of the restricted sheaf on the fiber $Y_b$.
\end{notation}

\section{Cycle spaces, Gauduchon forms, and integral classes}
\label{s:revised-preliminaries}

We organize the preliminaries around the componentwise proof outlined
in the introduction.  Section~\ref{ss:revised-cycle-tools} fixes the
cycle-space framework, the compactness criterion, and Barlet's
algebraic-dimension theorem.  Section~\ref{ss:revised-gauduchon-forms}
provides the forms and integration tools for divisor masses.
Section~\ref{ss:revised-integral-obstruction} develops the fixed-class
approach: base change, total-space line bundles, finite-order torsion,
and mass bounds where the obstruction vanishes.  The complementary
good-filling constructions are separated in
Section~\ref{ss:revised-good-fillings}.

\subsection{Divisor cycle spaces and componentwise compactness}
\label{ss:revised-cycle-tools}

We use Barlet spaces to retain both the supports and the positive
integral multiplicities of divisors; see \cite{BM19} for a systematic
treatment.  The tools below are used in
Section~\ref{ss:componentwise-punctured} to prepare a horizontal
component and in Section~\ref{ss:componentwise-conclusion} to conclude
properness and the algebraic-dimension inequality.

\subsubsection{Cycles and their integral classes}
Let $X$ be a compact complex $n$-dimensional manifold.  For
$0\le k\le n$, an \emph{effective $k$-cycle} is a finite formal sum
\[
 Z=\sum_{i\in I}m_iZ_i,
\]
where the $Z_i$ are pairwise distinct irreducible analytic subsets of
complex dimension $k$ and the $m_i$ are positive integers.  Its \emph{support}
is $|Z|=\bigcup_iZ_i$.  The \emph{Barlet space} $\cC_k(X)$ parametrizes these
cycles with their multiplicities, not merely their supports.

The integral homology class of $Z$ has a Poincar\'e-dual class
$c(Z)\in H^{2n-2k}(X,\Z)$.  If a smooth closed form $\eta_Z$
represents its real image, then
\[
 \int_Z\alpha=\int_X\alpha\wedge\eta_Z
\]
for every smooth closed $2k$-form $\alpha$.  For a divisor,
$c(Z)=c_1(\cO_X(Z))$.  We retain the integral class even though its
smooth representative detects only its real image.

For the family \eqref{eq:family}, Ehresmann's fibration theorem
\cite{Eh47} and contractibility of the disc give a differentiable
trivialization with common fiber $M$.  In particular, $R^2\pi_*\Z$ is
constant, so integral divisor classes on different fibers can be compared
in $H^2(M,\Z)$.  This is the \emph{Ehresmann trivialization} used
throughout the paper; it is not a holomorphic trivialization.

\subsubsection{Absolute and relative divisor spaces}
\begin{thm}[Compactness of divisor components,
{\cite[Theorem~6.0.13]{Ba15}}]
\label{thm:revised-absolute-divisors}
Every irreducible component of $\cC_{n-1}(X)$ is compact.
\end{thm}

Let $f:\cX\to B$ be a smooth family of connected compact complex
$n$-dimensional manifolds.  The \emph{relative Barlet space}
$\cC_k(\cX/B)$ parametrizes effective compact $k$-cycles contained
in the fibers, with the \emph{canonical projection}
\[
 \mu:\cC_k(\cX/B)\longrightarrow B.
\]
We regard its points as pairs $(b,Z)$, so that the zero cycle, if
included, retains its base point.  For the family \eqref{eq:family},
write $\mathscr C:=\cC_{n-1}(\cX/\Delta)$.
For a proper $n$-equidimensional map $f:Y\to B$, we similarly write
$\cC_{n-1}(f)=\cC_{n-1}(Y/B)$ for its \emph{relative divisor cycle space}.
An \emph{irreducible component} always means a maximal irreducible closed
analytic subset of this space.

\begin{thm}[Barlet classification and fundamental cycles,
{\cite[Theorems~1.8--1.9, Main Theorem, and Corollary~3.3]{Ma07}}]
\label{thm:revised-cycle-universal}
Let $M$ be a complex space.
\begin{enumerate}[$(i)$]
\item The reduced Barlet space $\cC_k(M)$ has a universal proper
analytic family of compact effective $k$-cycles.  A proper analytic
family over a reduced space $S$ determines a holomorphic map
$S\to\cC_k(M)$ sending $s$ to its cycle; conversely, pulling back the
universal family recovers that family, with multiplicities.
\item For a flat proper family of pure $k$-dimensional subspaces,
the fundamental cycles form an analytic family over the reduced base.
\end{enumerate}
\end{thm}

An analytic family of cycles retains the integer degrees of its local
multigraphs.  It need not be a flat family of analytic subspaces; the
last assertion of the theorem gives the passage from flat families to
cycle families that is used below.

\begin{lemma}[Constancy of the cycle class]
\label{lem:revised-cycle-class}
Under an Ehresmann trivialization, the integral homology class of an analytic
family of relative cycles is locally constant.  Hence every connected
analytic family of divisors has one Poincar\'e-dual class in $H^2(M,\Z)$.
\end{lemma}

\begin{proof}
For the definition and existence of analytic relative fundamental classes,
see \cite[\S~2~B, Definition~1, and \S~2~D, Theorem~2]{Ba80}.
The integral homology statement is obtained by the
following topological argument.
Along a path in the parameter space, adapted scales sweep an integral bordism
between the endpoint cycles; the constant local degrees retain the
multiplicities.  The endpoint cycles are therefore homologous in the
fixed differentiable fiber determined by the trivialization.
\end{proof}

\subsubsection{Compactness, proper images, and density}
Bishop compactness is used in its cycle-space form, which controls both
the supports and the positive integral multiplicities.

\begin{lemma}[Cycle compactness from a mass bound]\label{lem:revised-bishop-cycle}
Let $K\Subset\Delta$ be compact and let $H$ be a fixed Hermitian metric on
$\cX$.  A closed subset $A\subset\cC_{n-1}(\cX/\Delta)$ whose cycles are
supported in $\pi^{-1}(K)$ and have uniformly bounded $H$-volume is
compact.
\end{lemma}

\begin{proof}
This lemma is well known to experts; we include the proof for the
reader's convenience.
Properness of $\pi$ makes $\pi^{-1}(K)$ compact.  The compactness criterion
for Barlet spaces states that a family of compact $k$-cycles is relatively
compact when all supports lie in one compact subset and their volumes for
one Hermitian metric are uniformly bounded, by
\cite[Theorem~4.2.72]{BM19}.

That theorem is the cycle-theoretic form of Bishop compactness.  Its
adapted-scale formulation retains the positive integral multiplicities,
so its conclusion is convergence in the Barlet space rather than only
Hausdorff convergence of supports.

The assumptions of the criterion hold for $A$.  Thus $A$ is relatively
compact in $\cC_{n-1}(\cX/\Delta)$.  Since $A$ is closed and the Barlet
space is Hausdorff, $A$ is compact.
\end{proof}

Properness then gives analytic images by Remmert's theorem.

\begin{thm}[Proper mapping theorem,
{\cite[\S~7, Satz~1]{Gr60}}]\label{thm:revised-remmert}
Let $f:Y\to B$ be a proper holomorphic map of complex spaces.  For every
closed analytic subset $Z\subset Y$, its image $f(Z)$ is a closed
analytic subset of $B$.
\end{thm}

In particular, the image of a proper map to a connected Riemann surface
is the whole surface if it contains a nonempty open subset.
The following density statement applies before properness is known.

\begin{lemma}\label{lem:revised-open-density}
Let $Y$ be an irreducible reduced complex space and let
$h:Y\to B$ be a nonconstant holomorphic map to a Riemann surface.
Then $h$ is open.  If $V_B\subset B$ is dense and open, then
$h^{-1}(V_B)$ is dense in $Y$.
\end{lemma}

\begin{proof}
This follows directly from \cite[Chapter~II, Theorem~(5.7)]{Dem}.
For completeness, we give a proof.
Fix $y\in Y$.  No neighbourhood of $y$ can be mapped to a point: otherwise
$h$ would be constant on a nonempty open subset and hence, by the identity
principle on the irreducible reduced space, constant on $Y$.  Consequently
$h$ cannot be constant on every local irreducible branch through $y$, for
then it would be constant on a neighbourhood of $y$.  Choose a branch on
which its germ is nonconstant.  A generic linear section through $y$ gives
a complex curve germ $C\subset Y$ through $y$ on which $h$ is nonconstant.
After normalization, the ordinary open mapping theorem for one complex variable
shows that $h(C\cap W)$ contains a neighbourhood of $h(y)$ for every
sufficiently small neighbourhood $W$ of $y$.  Hence $h$ is open.
If $h^{-1}(V_B)$ were not dense, some nonempty open $W\subset Y$ would map
into $B\setminus V_B$.  Openness of $h$ would give a nonempty open subset
of $B$ contained in the nowhere-dense complement, a contradiction.
\end{proof}

\subsubsection{From componentwise properness to algebraic dimension}
The following theorem is the final input in the proof of
Theorem~\ref{thm:main}.  It reduces the algebraic-dimension conclusion
to properness of the relative divisor components, without requiring a
separate choice or completion of a good filling.

\begin{thm}[Local algebraic-dimension theorem,
{\cite[Proposition~4.0.8]{Ba15}}]
\label{thm:revised-barlet-local}
Let $f:Y\to B$ be a proper surjective $n$-equidimensional holomorphic map
between irreducible complex spaces.  Assume that every irreducible
component of $\cC_{n-1}(f)$ is proper over $B$.  Then there are a
nonnegative integer $a$ and a countable union $\Sigma$ of closed
irreducible analytic subsets of $B$, each with empty interior, such that
\[
 a(Y_b)=a\quad(b\in B\setminus\Sigma),
 \qquad a(Y_b)\ge a\quad(b\in B).
\]
\end{thm}

\subsection{The \texorpdfstring{$\partial\bar\partial$}{ddbar}-lemma and Gauduchon forms}
\label{ss:revised-ddbar-locus}
\label{ss:revised-gauduchon-forms}

This subsection supplies the forms and integration identities used
in Section~\ref{ss:componentwise-mass}.  The current
$\partial\bar\partial$-lemma makes it possible to compare a divisor with
a smooth representative of its class.  A smooth Gauduchon family then
annihilates the resulting exact correction, while continuity of cycle
integration transfers mass bounds to limiting cycles.

\subsubsection{The smooth and current formulations}
\begin{defn}\label{def:revised-ddbar}
A compact complex manifold $X$ satisfies the
\emph{$\partial\bar\partial$-lemma} if every $d$-closed pure-type smooth form which
is $d$-exact is $\partial\bar\partial$-exact.  Equivalently, for every
bidegree,
\[
 \ker d\cap\operatorname{im}d
 =\operatorname{im}\partial\bar\partial
\]
inside the space of pure-type forms.  The equivalent current formulation
will also be used.
\end{defn}

The following current form is the precise consequence needed in the mass
calculation.

\begin{lemma}[{\cite[Lemma~5.0.12]{Ba15}}]\label{lem:revised-current-ddbar}
Let $X$ be a compact connected complex manifold satisfying the
$\partial\bar\partial$-lemma.  If $T$ is a real, $d$-exact current of type
$(1,1)$ on $X$, then there is a real $0$-current $\Phi$ such that
\begin{equation}\label{eq:revised-current-ddbar}
 T=i\partial\bar\partial\Phi.
\end{equation}
\end{lemma}

\begin{proof}
Barlet proves that any $d$-exact $(1,1)$-current can be written
$T=\partial\bar\partial\Theta$ for a complex $0$-current $\Theta$.  When $T$ is real, conjugation gives
$-\partial\bar\partial\overline\Theta=T$.  Hence
$\Phi=(\Theta-\overline\Theta)/(2i)$ is real and satisfies
\eqref{eq:revised-current-ddbar}.
\end{proof}

\subsubsection{Gauduchon forms and continuity of the mass}
We choose a smooth Gauduchon family on the entire disc, including the
possibly exceptional fiber.  The strongly Gauduchon condition is recalled
as well for the complementary constructions in
Section~\ref{ss:revised-good-fillings}.

\begin{defn}\label{def:revised-gauduchon}
Let $X$ be a compact complex $n$-fold.  A Hermitian form $\gamma$ is
\emph{Gauduchon} if
\[
 \partial\bar\partial\gamma^{n-1}=0.
\]
It is \emph{strongly Gauduchon} if
$\partial\gamma^{n-1}$ is $\bar\partial$-exact.  Equivalently, there is a
real $d$-closed $(2n-2)$-form $\Omega$ such that
$\Omega^{n-1,n-1}=\gamma^{n-1}$, by
\cite[Definition~4.1 and Proposition~4.2]{Po13}.
\end{defn}

Every $\partial\bar\partial$-manifold is strongly Gauduchon: start with a
Gauduchon metric and apply the $\partial\bar\partial$-lemma to the relevant
$\partial$-exact form.  This implication is also recorded in
\cite[p.~530]{Po13}.

\begin{lemma}[{\cite{Ga77}; \cite[\S~3]{Po13}}]
\label{lem:revised-smooth-gauduchon}
There is a smooth family of positive $(n-1,n-1)$-forms
\begin{equation*}
 G_t=\gamma_t^{n-1},\qquad
 \partial_t\bar\partial_tG_t=0,
 \qquad t\in\Delta.
\end{equation*}
\end{lemma}

\begin{lemma}[{\cite[Proposition~4.2.17.(2)]{BM19}}]
\label{lem:revised-absolute-G}
After choosing a smooth horizontal splitting of
$T\cX\to\pi^*T\Delta$, the relative form $G=(G_t)$ of
Lemma~\ref{lem:revised-smooth-gauduchon} extends to a smooth absolute real
$(2n-2)$-form $\widetilde G$ on $\cX$ satisfying
$\widetilde G|_{X_t}=G_t$.  In particular,
\begin{equation*}
 \cC_{n-1}(\cX/\Delta)\longrightarrow\R,
 \qquad Z\longmapsto\int_Z\widetilde G=\int_ZG_{\pi(Z)}
\end{equation*}
is continuous.
\end{lemma}

\subsection{The integral-class obstruction and mass bounds}
\label{ss:revised-integral-obstruction}

Following the fixed-integral-class approach discussed in the
introduction, we encode the existence of a total-space line bundle by
a holomorphic obstruction section.  Base change gives its initial
vanishing, and the exponential sequence gives line bundles where it
vanishes.  The crucial remaining information is finite-order
annihilation, not necessarily global vanishing; this is the input for
the phase argument in Section~\ref{ss:componentwise-mass}.

\subsubsection{Coherence and base change on the
\texorpdfstring{$\partial\bar\partial$}{ddbar}-locus}
We first record several basic results from Grauert's direct image theory used in
Section~\ref{ss:componentwise-punctured} to pass from fiberwise
vanishing to vanishing of a holomorphic obstruction section.

\begin{thm}[Coherence, semicontinuity, and base change,
{\cite[Hauptsatz~I and \S~7, S\"atze~3 and~5]{Gr60}},
{\cite[\S~10.5.4]{GR84}}]
\label{thm:revised-grauert-tools}
Let $f:Y\to B$ be a proper holomorphic map and let $\cF$ be a coherent
analytic sheaf on $Y$.
\begin{enumerate}[$(i)$]
\item \textup{(Coherence)} Every $R^qf_*\cF$ is coherent.
\item \textup{(Semicontinuity)} If $B$ is a complex manifold and $\cF$
is flat over $B$, the function
$b\mapsto h^q(Y_b,\cF|_{Y_b})$ is upper semicontinuous in the analytic
Zariski topology: for every integer $d\ge0$, the locus where
$h^q(Y_b,\cF|_{Y_b})\ge d$ is a closed analytic subset of $B$.
\item \textup{(Base change)} Under the hypotheses of \textup{(ii)}, where
the function $b\mapsto h^q(Y_b,\cF|_{Y_b})$ is locally constant,
$R^qf_*\cF$ is locally free and the natural
base-change map is an isomorphism:
\[
 (R^qf_*\cF)_b/\mathfrak m_b(R^qf_*\cF)_b
 \longrightarrow H^q(Y_b,\cF|_{Y_b}).
\]
Here $\mathfrak m_b$ is the maximal ideal of $\cO_{B,b}$.
\end{enumerate}
\end{thm}

In our applications, $f$ is a
submersion and $\cF$ is a vector bundle, so flatness holds.  Coherence
alone does not identify a stalk modulo $\mathfrak m_b$ with fiber
cohomology: that identification is the additional base-change assertion.

\begin{lemma}[Base change over a disc of
$\partial\bar\partial$-fibers]\label{lem:revised-ddbar-base-change}
Let $U\subset\Delta$ be a complex disc such that every $X_t$, $t\in U$,
satisfies the $\partial\bar\partial$-lemma.  Then, for all integers
$0\le p,q\le n$, the sheaves
$R^q\pi_*\Omega^p_{\cX/\Delta}|_U$ are locally free and satisfy base change.
\end{lemma}

\begin{proof}
For every integer $0\le k\le 2n$, the $\partial\bar\partial$-lemma gives
\begin{equation*}
 b_k(M)=\sum_{\substack{p+q=k\\0\le p,q\le n}}
 h^q(X_t,\Omega^p_{X_t})\qquad(t\in U).
\end{equation*}
Each summand is upper semicontinuous, while the sum is constant.  Shrinking
about any point of $U$ shows that every summand is locally constant.
Theorem~\ref{thm:revised-grauert-tools} then makes $R^q\pi_*\Omega^p_{\cX/\Delta}|_U$ locally free
and gives base change, as in
\cite[Satz~5]{Gr60}.
\end{proof}

\subsubsection{The obstruction section and total-space line bundles}
\begin{thm}[Cartan's theorem B,
{\cite[Theorem~B]{Ca53}}]\label{thm:revised-cartan}
Let $W$ be a Stein manifold and let $\cF$ be a coherent analytic sheaf
on $W$.  Then $H^q(W,\cF)=0$ for every $q\ge1$.
\end{thm}

Discs and punctured discs are Stein.  The theorem removes higher
cohomology on the base without requiring the total space $\cX_W$
to be Stein.

For clarity, the Leray spectral sequence for a proper map $f:Y\to W$
and a coherent sheaf $\cF$ starts with
\[
 E_2^{p,q}=H^p(W,R^qf_*\cF)\ \Longrightarrow\ H^{p+q}(Y,\cF).
\]
Theorem~\ref{thm:revised-grauert-tools} makes the direct images coherent, and
Theorem~\ref{thm:revised-cartan} ensures that all columns with $p>0$ vanish.  Hence
\[
 H^q(Y,\cF)\simeq H^0(W,R^qf_*\cF).
\]
This is the precise cohomological passage used in the lifting argument
below, not an assertion that fiberwise vanishing alone implies vanishing
of a section of a coherent sheaf.

\begin{defn}[Torsion and torsion-freeness, {\cite[\S~7.5]{Re94}}]
\label{def:revised-torsion-sheaf}
Let $\cF$ be a coherent sheaf on a disc $D$.  Its \emph{torsion
subsheaf} $T(\cF)$ has stalks
\begin{equation*}
 T(\cF)_p:=\{s\in\cF_p:\text{there is }0\ne g\in\cO_{D,p}
                         \text{ with }gs=0\}.
\end{equation*}
The sheaf $\cF$ is \emph{torsion-free at $p$} if $T(\cF)_p=0$,
and \emph{torsion-free} if this holds at every point.
\end{defn}

The finite-order vanishing needed below is established directly in
Lemma~\ref{lem:revised-discrete-support}.

A class $c\in H^2(M,\Z)$ therefore determines a total
topological class $\widetilde c\in H^2(\cX,\Z)$ and a holomorphic section
\begin{equation*}
 o(c)\in H^0\!\left(\Delta,R^2\pi_*\cO_{\cX}\right)
\end{equation*}
through the morphism $\Z\to\cO_{\cX}$.

For an open disc $W\subset\Delta$ and $t\in W$, restriction and the natural Leray
isomorphism give the following commutative diagram, extending
\cite[equation~(4.6)]{RT21}:
\begin{equation}\label{eq:revised-exponential-restriction}
\begin{CD}
 \Pic(\cX_W) @>{c_1}>> H^2(\cX_W,\Z)
 @>{\iota_W}>> H^2(\cX_W,\cO_{\cX_W})
 @>{\lambda_W}>{\sim}> (R^2\pi_*\cO_{\cX})(W)\\
 @VV{\mathrm{res}_t}V @VV{\mathrm{res}_t}V @VV{\mathrm{res}_t}V
 @VV{\rho_t}V\\
 \Pic(X_t) @>{c_1}>> H^2(X_t,\Z)
 @>{\iota_t}>> H^2(X_t,\cO_{X_t}) @= H^2(X_t,\cO_{X_t}).
\end{CD}
\end{equation}
Here $\iota_W$ and $\iota_t$ are induced by $\Z\to\cO$.
The first three terms of each row come from the exponential sequence
$0\to\Z\to\cO\to\cO^*\to1$, and are exact at the integral
cohomology term.  The first three vertical arrows are restriction maps;
naturality of the exponential sequences gives the first two squares.
Coherence and Cartan's theorem~B on the Stein disc $W$ give the
Leray isomorphism $\lambda_W$, where
$(R^2\pi_*\cO_{\cX})(W)=H^0(W,R^2\pi_*\cO_{\cX})$.
The last vertical arrow is
\[
 \begin{aligned}
 \rho_t:(R^2\pi_*\cO_{\cX})(W)
 &\longrightarrow (R^2\pi_*\cO_{\cX})_t\\
 &\longrightarrow
 (R^2\pi_*\cO_{\cX})_t/\mathfrak m_t(R^2\pi_*\cO_{\cX})_t
 \xrightarrow{\mathrm{bc}_t}H^2(X_t,\cO_{X_t}),
 \end{aligned}
\]
where $\mathfrak m_t\subset\cO_{\Delta,t}$ is the maximal ideal and
$\mathrm{bc}_t$ is the natural base-change map.  Naturality gives
$\rho_t\circ\lambda_W=\mathrm{res}_t$ on total-space cohomology.
Neither torsion-freeness of
$R^2\pi_*\cO_{\cX}$ nor invertibility of $\mathrm{bc}_t$ is required.
The Ehresmann trivialization identifies $\widetilde c|_{\cX_W}$ with $c$ under the
second vertical arrow.  The corresponding images are therefore
\[
 \lambda_W\!\left(\iota_W(\widetilde c|_{\cX_W})\right)=o(c)|_W,
 \qquad
 \rho_t(o(c)|_W)=\iota_t(c).
\]

Vanishing of $\rho_t(o(c)|_W)$ need not imply vanishing of the germ
$o(c)_t$: for example, $z\bmod z^2$ is nonzero in $R/(z^2)$ but
has zero image in $R/(z)$, where $R=\cO_{\Delta,0}$.
Thus fiberwise vanishing alone does not justify the total-space lift;
the following lemma requires vanishing of the section on $W$.

\begin{lemma}\label{lem:revised-local-lift}
Let $W\subset\Delta$ be an open disc, possibly $W=\Delta$.
If $o(c)|_W=0$, then there is a line
bundle $\cL_W\in\Pic(\cX_W)$ with
$c_1(\cL_W)=\widetilde c|_{\cX_W}$.  A Hermitian metric on $\cL_W$ gives a
smooth real $d$-closed $(1,1)$-form $\eta$ on $\cX_W$ whose restriction
$\eta_t$ represents the real image of $c$.
\end{lemma}

\begin{proof}
Theorems~\ref{thm:revised-grauert-tools} and~\ref{thm:revised-cartan}
make the Leray spectral sequence
over the Stein disc $W$ degenerate for $\cO$ in the form
\[
 H^2(\cX_W,\cO_{\cX_W})
 \simeq H^0(W,R^2\pi_*\cO_{\cX}).
\]
Thus the image of $\widetilde c$ in $H^2(\cX_W,\cO)$ is zero.  Exactness
of the exponential sequence
\[
 0\longrightarrow\Z\longrightarrow\cO
 \xrightarrow{\exp(2\pi i\,\cdot)}\cO^*\longrightarrow1
\]
produces $\cL_W$, as in the top row of
\eqref{eq:revised-exponential-restriction}.  The normalized Chern curvature of a
smooth Hermitian metric has the asserted properties.
\end{proof}

\subsubsection{Discrete support and finite-order torsion}
The next lemma is the key sheaf-theoretic input for
Sections~\ref{ss:componentwise-punctured} and~\ref{ss:componentwise-mass}:
it controls the obstruction even when the direct-image sheaf is not
torsion-free.

Set $\cF:=R^2\pi_*\cO_{\cX}$.  The notation $\cO_\Delta o(c)$ denotes
the \emph{cyclic subsheaf} of $\cF$ generated by $o(c)$, namely
\[
 \cO_\Delta o(c):=
 \operatorname{im}\!\left(\cO_\Delta
 \xrightarrow{\,f\mapsto f o(c)\,}\cF\right).
\]
Its stalk at $p$ is $\cO_{\Delta,p}o(c)_p$.
\begin{lemma}[Discrete support and finite torsion order;
cf.\ {\cite[Proposition~4.5]{RT21}}]\label{lem:revised-discrete-support}
If $o(c)$ vanishes on a nonempty open subset of $\Delta$, then
\begin{equation*}
 S(c):=\Supp\bigl(\cO_\Delta\,o(c)\bigr)
\end{equation*}
is a possibly empty discrete complex-analytic subset.  For every
$p\in\Delta$, with local coordinate $z=t-p$, there is a least integer
$N(p)\ge0$ such that
\begin{equation}\label{eq:revised-torsion-order}
 z^{N(p)} o(c)=0
\end{equation}
near $p$.  Here $N(p)=0$ if and only if $p\notin S(c)$; at a point of
$S(c)$ one necessarily has $N(p)\ge1$.  In particular, $o(c)=0$ on
$\Delta$ if and only if $S(c)=\varnothing$, equivalently if $N(p)=0$
for every $p$.
\end{lemma}

\begin{proof}
$\cF$ is coherent by Theorem~\ref{thm:revised-grauert-tools}, so its cyclic
subsheaf $\cN:=\cO_\Delta o(c)$ is coherent.  Its support is the zero set
of its annihilator ideal sheaf and is therefore analytic.  It is proper
because $o(c)$ vanishes on a nonempty open set.  Since $\Delta$ is a
connected complex curve, this proper analytic subset is discrete.

If $p\notin S(c)$, the germ $o(c)_p$ is zero, hence the section vanishes
on a neighbourhood of $p$.  Thus the least exponent there is $N(p)=0$.

Fix $p\in S(c)$ and put $R:=\cO_{\Delta,p}\simeq\C\{z\}$, where
$z=t-p$.  The ring $\C\{z\}$ consists of convergent power-series germs.
Every nonzero germ has the form $z^m u$, where $m\ge0$ is its vanishing
order and $u$ is a unit.  Choosing the least vanishing order among the
nonzero elements of an ideal shows that every nonzero ideal is $(z^m)$.
This is the discrete valuation ring property used here.

The cyclic stalk has the presentation
\[
 \cN_p\simeq R/I_p,\qquad
 I_p:=\{g\in R:g\,o(c)_p=0\}.
\]
Here $I_p\ne R$ because $\cN_p\ne0$.  Also $I_p\ne0$: otherwise the
coherent kernel of $\cO_\Delta\to\cN$ would vanish near $p$, giving
$\cN\simeq\cO_\Delta$ there, contrary to the isolated support.  Hence
$I_p=(z^N)$ for some $N\ge1$, and
\[
 \cN_p\simeq R/(z^N),\qquad z^N\cN_p=0.
\]
In particular, $\cN_p$ is a finitely generated torsion module, of length
$N$.  The equality $z^No(c)_p=0$ is an equality of germs, so
$z^No(c)=0$ on a sufficiently small neighbourhood of $p$.
Since $I_p=(z^N)$, no smaller power annihilates the generator, so
$N=N(p)$.  This proves \eqref{eq:revised-torsion-order} and the assertions about
the zero exponent.
\end{proof}

Rao--Tsai's supported-section criterion
\cite[Proposition~4.5]{RT21} uses the coherent-sheaf support theorem
\cite[Chapter~II, Proposition~3.2]{BS76} and local freeness on the
punctured base.  Here only finite-order annihilation is needed:
the cyclic-subsheaf argument above proves it without local freeness
of the ambient sheaf away from $S(c)$.

The zero-order case deserves separate emphasis.  Initial vanishing on
a nonempty open set does not by itself give $o(c)=0$ on all of $\Delta$:
the remaining germs may be nonzero torsion.  If $\cF$ is torsion-free,
however, \eqref{eq:revised-torsion-order} forces every germ to vanish.  Then
$S(c)=\varnothing$ and $N(p)=0$ everywhere.  Conversely, a particular
$o(c)$ can vanish even when $\cF$ has nonzero torsion.

In Rao--Tsai's argument \cite[Propositions~4.12 and~4.13]{RT21},
constancy of $h^{0,1}(X_t)$ supplies torsion-freeness, so initial
vanishing extends globally.  This is the case $N(p)=0$ for all $p$;
Lemma~\ref{lem:revised-local-lift} with $W=\Delta$ then gives the
total-space line bundle without a positive-order phase argument.
\subsubsection{Mass bounds away from the obstruction}
We now combine the integral-class lifting criterion with the Gauduchon
forms of Section~\ref{ss:revised-gauduchon-forms}.  The resulting mass
identity is used in Section~\ref{ss:componentwise-punctured} and applies
to any irreducible component, not only one obtained from a good filling.

\begin{lemma}\label{lem:revised-mass-off-S}
Let $\Lambda\subset\mathscr C$ be an irreducible component with fixed
integral divisor class $c$, and suppose $o(c)|_D=0$ on a disc
$D\subset\Delta$.  Let $\eta$ be the curvature representative given
by Lemma~\ref{lem:revised-local-lift}.  Whenever $t\in D$ and $X_t$ satisfies
the $\partial\bar\partial$-lemma, every $Z\in\Lambda_t$ satisfies
\begin{equation}\label{eq:revised-mass-off-S}
 \int_ZG_t=\int_{X_t}\eta_t\wedge G_t.
\end{equation}
The right-hand side is smooth in $t$, so these masses are uniformly
bounded on compact subsets of $D$, independently of $Z$.
\end{lemma}

\begin{proof}
The currents $[Z]$ and $\eta_t$ are real, closed, of type $(1,1)$,
and represent the same real class.  The current
$\partial\bar\partial$-Lemma~\ref{lem:revised-current-ddbar} gives
$[Z]-\eta_t=i\partial_t\bar\partial_t\Phi_Z$ for a real $0$-current
$\Phi_Z$.  Pairing with $G_t$ and using
$\partial_t\bar\partial_tG_t=0$ gives \eqref{eq:revised-mass-off-S}.
Smooth fiber integration proves the last assertion.
\end{proof}

If $\mu|_\Lambda$ is nonconstant and the
$\partial\bar\partial$-parameters form a dense open subset of the
base, Lemma~\ref{lem:revised-open-density} makes their inverse image dense in
$\Lambda$.  Continuity of cycle integration
(Lemma~\ref{lem:revised-absolute-G}) extends each local mass bound to all
cycles over a smaller disc.  Positivity of $G_t$ compares these masses
with Hermitian volumes on compact sets, and
Lemma~\ref{lem:revised-bishop-cycle} gives properness over any open region
where $o(c)$ vanishes.  The possible remaining points of $S(c)$ are
handled by the phase estimate in
Section~\ref{s:componentwise-main-proof}.

\subsection{Strongly Gauduchon properness and good fillings}
\label{ss:revised-good-fillings}

The main argument in Section~\ref{s:componentwise-main-proof} uses the
componentwise criterion of Theorem~\ref{thm:revised-barlet-local}
directly.  We collect here the complementary geometric tools used in
Remark~\ref{rem:componentwise-good-filling-route}, the examples, and
Appendix~\ref{s:applications}: strongly Gauduchon properness, relative
Kodaira images, and continuation of a local filling to a global component.

D. Lieberman \cite{Lb78} appears to have been the first to notice the ``compactness of the Chow scheme'' in the K\"ahler context
and its strong implications for K\"ahler geometry, initially for the structure of automorphism groups. These observations were the starting point for later developments\footnote{We thank F. Campana for pointing this out to us.}.

F. Campana \cite{Ca81} studied the algebraic reduction of a proper weakly K\"ahler morphism $\phi$ via the relative Kodaira map, proving the semicontinuity of the algebraic dimension of the fibers of $\phi$. Notice that the weakly K\"ahler assumption implies that all fibers of $\phi$ are in Fujiki class $\mathcal{C}$ and that the
irreducible components of the relative cycle spaces are proper for cycles
of any dimension \cite[Theorem 4.5 and Proposition 4.8]{Fj78}, while Barlet \cite{Ba15} used only properness for relative codimension-one cycles and did not impose a K\"ahler-type assumption on the fibers.

\subsubsection{Strongly Gauduchon properness}
A \emph{relative strongly Gauduchon form} is a smooth real relative
$(2n-2)$-form whose restriction to each fiber is closed and has a strictly
positive $(n-1,n-1)$ component.  A strongly Gauduchon form on one fiber
extends to such a relative form after shrinking the base, by
\cite[Lemma~3.0.7 and the following remark]{Ba15}.  If every fiber is
strongly Gauduchon, these local extensions can be patched by a partition
of unity on the base.

\begin{lemma}[Strongly Gauduchon properness,
{\cite[Proposition~3.0.6]{Ba15}}]
\label{prop:revised-barlet-proper}
If $f:\cX\to B$ admits a relative strongly Gauduchon form, then
the restriction of $\mu:\cC_{n-1}(\cX/B)\to B$ to each connected
component, and hence to each irreducible component, is proper over $B$.
\end{lemma}

\subsubsection{Good fillings and algebraic dimension}
Let $\Gamma$ be an irreducible $B$-proper analytic subset of
$\cC_{n-1}(\cX/B)$.  Its \emph{tautological incidence graph} is
\[
 G_\Gamma:=\{(Z,x)\in\Gamma\times_B\cX:x\in|Z|\}.
\]

\begin{defn}[Good filling, {\cite[\S~2.2]{Ba15}}]
\label{def:revised-good-filling}
The family $\Gamma$ is a \emph{good filling} if:
\begin{enumerate}[$(i)$]
\item the evaluation $G_\Gamma\to\cX$ is surjective;
\item the generic member of the tautological family is irreducible.
\end{enumerate}
A good filling need not be an irreducible component of the cycle space.
Conversely, being a component does not imply either of these conditions
or properness over the base.
\end{defn}

For a good filling $\Gamma$, after a proper $B$-modification
$\alpha_\Gamma:\cX_\Gamma\to\cX$, Barlet constructs a holomorphic
relative Kodaira map
\[
 K_{\Gamma/B}:\cX_\Gamma\longrightarrow
 \cC_{\dim\Gamma-\dim B-1}(\Gamma/B).
\]
It is obtained from the fiber map of the incidence graph followed by
the direct-image map to $\Gamma$; at a general point $x$, it records
the incidence cycle of divisors containing $x$.
Its image $Q_\Gamma$ is proper over $B$ by
\cite[Lemma~2.2.2]{Ba15}.  Equality of its generic fiber dimension with
the generic algebraic dimension requires a filling chosen to give an
algebraic reduction.  For an arbitrary good filling we use only the
following inequality.

\begin{prop}[Good-filling inequality]
\label{prop:revised-barlet-endgame}
Let $B$ be a connected complex curve and let
$\Gamma\subset\cC_{n-1}(\cX/B)$ be a $B$-proper good filling.
If the generic fiber of $Q_\Gamma\to B$ has dimension at least $a$,
then $a(X_b)\ge a$ for every $b\in B$.
\end{prop}

\begin{proof}
Barlet proves in \cite[Lemma~2.2.2]{Ba15} that $Q_\Gamma\to B$ is proper.
His exceptional set $T_\Gamma$ has
codimension at least two and is therefore empty when $\dim B=1$, by
\cite[Lemma~2.2.1]{Ba15}.  Let $q$ be the generic fiber
dimension of the proper map $Q_\Gamma\to B$.  On a dense open subset its
fibers have dimension $q$, and $q\ge a$ by hypothesis.  Apply Barlet's
fiber inequality \cite[Lemma~2.2.3]{Ba15} to that dense set;
his conclusion
is $a(X_b)\ge q\ge a$ for every $b\in B$.
\end{proof}

\subsubsection{Local fillings and continuation of components}
For the alternative good-filling argument, assume that every fiber
over $\Delta^*$ satisfies the $\partial\bar\partial$-lemma and that
\[
 a_*:=\inf_{t\in\Delta^*}a(X_t)>0.
\]
Choose a relatively compact disc $U\Subset\Delta^*$.

\begin{prop}[The local good filling]\label{prop:revised-local-filling}
After shrinking $U$, there are an integer $a_U\ge a_*$ and an irreducible component
\begin{equation*}
 \Gamma_U\subset\cC_{n-1}(\cX_U/U)
\end{equation*}
such that:
\begin{enumerate}[$(i)$]
\item $\Gamma_U\to U$ is proper and $\Gamma_U$ is a good filling;
\item the generic fiber of the relative Kodaira image
$Q_{\Gamma_U}\to U$ has dimension $a_U$;
\item $a(X_t)=a_U$ away from a countable union of proper analytic subsets
of $U$, and $a(X_t)\ge a_U$ for every $t\in U$.
\end{enumerate}
\end{prop}

\begin{proof}
Every fiber over $U$ is strongly Gauduchon.  A relative strongly
Gauduchon form therefore exists by
\cite[Lemma~3.0.7 and the following remark]{Ba15}.
Lemma~\ref{prop:revised-barlet-proper} makes every connected component of the
relative divisor space proper over $U$.

Every irreducible component, being closed in one of these connected
components, is therefore proper.  We may therefore apply
\cite[Proposition~4.0.8 and its proof]{Ba15}.

By the Baire theorem, choose $t_1\in U$ outside the countable union of
exceptional loci in that proof, including the bigness and non-dominance
loci.  Choose the fiberwise good filling
$\Gamma_1\subset\cC_{n-1}(X_{t_1})$ which gives the algebraic reduction
there, and let $\Gamma_U$ be the relative irreducible component selected
by Barlet.

This component dominates $U$, has generically irreducible members, and its
relative Kodaira image has generic dimension
\[
 a_U=a(X_{t_1}).
\]

We check the first condition in Definition~\ref{def:revised-good-filling}, since
it is not a consequence of the word ``component.''  Barlet observes in
\cite[\S~2.2]{Ba15} that
the tautological incidence graph of an irreducible family whose generic
member is irreducible is itself irreducible.  The evaluation of that graph in $\cX_U$ is
proper; hence its image $Y$ is an irreducible closed analytic subset by
Theorem~\ref{thm:revised-remmert}.  It dominates $U$ and contains the whole fiber
$X_{t_1}$ because $\Gamma_1\subset(\Gamma_U)_{t_1}$ is a fiberwise good
filling.

Every fiber of the incidence graph contains an $(n-1)$-cycle, so
$\dim Y\ge n$.  If $\dim Y=n$, the inclusion $X_{t_1}\subset Y$ between
irreducible $n$-dimensional analytic sets would give $Y=X_{t_1}$,
contradicting dominance over $U$.  Hence $\dim Y=n+1$ and, since $\cX_U$ is
irreducible of that dimension, $Y=\cX_U$.  Thus $\Gamma_U$ is a good
filling.

Finally, Theorem~\ref{thm:revised-barlet-local} gives the assertions in $(iii)$.
Since $a(X_{t_1})\ge a_*$, we have $a_U\ge a_*$.
\end{proof}

The local family can be placed in a component over the original base
without assuming properness there.

\begin{lemma}\label{lem:revised-global-component}
Let $Y$ be a reduced complex space, let $Y^\circ\subset Y$ be open, and let
$A^\circ$ be an irreducible component of $Y^\circ$.  Then there is an irreducible component $A$ of $Y$ such that
\begin{equation}\label{eq:revised-component-restriction}
 A^\circ\subset A\cap Y^\circ.
\end{equation}
Moreover, $A^\circ$ is an irreducible component of $A\cap Y^\circ$,
$\dim A=\dim A^\circ$, and $A^\circ$ contains a nonempty subset which is
open in $A$.  In particular, with
\[
 \mathscr C:=\cC_{n-1}(\cX/\Delta),\qquad
 \mu:\mathscr C\to\Delta,
\]
there is an irreducible component $\Gamma\subset\mathscr C$ satisfying
\begin{equation*}
 \Gamma_U\subset\Gamma\cap\mu^{-1}(U),
\end{equation*}
and $\Gamma_U$ contains a nonempty open subset of $\Gamma$.
\end{lemma}

\begin{proof}
Complex spaces are locally Noetherian, and their local irreducible
components are locally finite.  In a neighbourhood meeting only finitely
many local components, the intersections of the local branches of
$A^\circ$ with the other components are proper analytic subsets.  We may
therefore choose a point $a\in A^\circ$ which lies on no other local
component of $Y^\circ$ and at which the germ of $A^\circ$ is irreducible.

A local irreducible component through $a$ is
contained in an irreducible component $A$ of $Y$.  Hence
$A\cap A^\circ$ is a closed analytic subset of the irreducible space
$A^\circ$ which contains a nonempty open subset.  The identity principle
gives $A^\circ\subset A$.

If an irreducible closed analytic subset of
$A\cap Y^\circ$ contains $A^\circ$, then it is also an irreducible closed
analytic subset of $Y^\circ$; maximality of $A^\circ$ forces equality.
Thus $A^\circ$ is an irreducible component of $A\cap Y^\circ$.

The other irreducible components of $A\cap Y^\circ$ form a locally finite
family, and their intersections with $A^\circ$ are proper analytic
subsets.  Choose a point of $A^\circ$ outside their union and shrink a
neighbourhood in $A\cap Y^\circ$.  The resulting nonempty subset is open
in $A$ and is contained in $A^\circ$.  In particular,
$\dim A=\dim A^\circ$.

One cannot in general strengthen \eqref{eq:revised-component-restriction} to
$A\cap Y^\circ=A^\circ$: for the irreducible curve $t=y^2$, the inverse
image of a small simply connected disc away from $t=0$ has two components.
Apply the valid containment statement to $Y=\mathscr C$ and
$Y^\circ=\mu^{-1}(U)$.
\end{proof}

\section{Phase compactness and componentwise properness}
\label{s:componentwise-main-proof}

We prove Theorem~\ref{thm:main} by establishing properness of every
irreducible component of the relative divisor space.  The integral-class
approach builds on Rao--Tsai's torsion-freeness method, which extends local vanishing of
the holomorphic obstruction to obtain a total-space line bundle when
$R^2\pi_*\cO_{\cX}$ is torsion-free.  Here a scalar phase estimate handles
the remaining finite-order torsion without that assumption and yields
a uniform bound for divisor masses.  Bishop compactness then gives
componentwise properness, and Barlet's Theorem~\ref{thm:revised-barlet-local}
yields the algebraic-dimension conclusion.  No good-filling assumption
is imposed on the component used in the mass estimate.

Throughout this section assume the hypotheses of Theorem~\ref{thm:main}.
An automorphism sends the possibly exceptional parameter to $0$.
We keep the original base $\Delta$; only auxiliary neighbourhoods are
shrunk.  Put $\mathscr C:=\cC_{n-1}(\cX/\Delta)$, let
$\mu:\mathscr C\to\Delta$ be the projection, and fix the smooth
Gauduchon family $G_t$ of Lemma~\ref{lem:revised-smooth-gauduchon}.
For an irreducible component $\Lambda$ of $\mathscr C$, endowed
with its reduced structure, write $\Lambda_t:=\Lambda\cap\mu^{-1}(t)$.
We call $\Lambda$
\emph{vertical} when $\mu|_\Lambda$ is constant and
\emph{horizontal} otherwise.  A discrete image of the connected space
$\Lambda$ is a single point.  Constants and auxiliary discs may depend
on $\Lambda$; no uniformity across different components is required.

\subsection{A scalar phase estimate}
\label{ss:componentwise-phase}

We isolate the analytic estimate used to control the possible singular
term in the divisor mass formula.  A lower bound for the real part of
$z^{-N}A(z)$ in every angular direction forces the whole expression to
remain bounded.  In the geometric application, this lower bound will
come from the nonnegativity of divisor masses.

\begin{lemma}[Phase lemma]\label{lem:componentwise-phase}
Let $N\ge0$ be an integer and let $A\in\mathcal C^N(D_0,\C)$,
where $D_0$ is a neighbourhood of $0$ in $\C$.  Suppose that there
are a real constant $C$ and $\varepsilon>0$ such that
\begin{equation}\label{eq:componentwise-phase-hypothesis}
 \operatorname{Re}\!\left(z^{-N}A(z)\right)\ge-C
 \qquad(0<|z|<\varepsilon).
\end{equation}
Then $A(z)=O(|z|^N)$, or equivalently $z^{-N}A(z)$ is bounded on
a punctured neighbourhood of $0$.
\end{lemma}

\begin{proof}
For $N=0$, continuity suffices.  Assume $N\ge1$ and write $z=x+iy$.
For a multi-index $\alpha=(\alpha_1,\alpha_2)$, use the notation
$|\alpha|=\alpha_1+\alpha_2$, $\alpha!=\alpha_1!\alpha_2!$,
$(x,y)^\alpha=x^{\alpha_1}y^{\alpha_2}$, and
$D^\alpha=\partial_x^{\alpha_1}\partial_y^{\alpha_2}$.
For $1\le m\le N$, apply the one-variable Taylor formula with integral
remainder \cite[formula~(1)]{KL03} to $g(s)=A(sx,sy)$, $0\le s\le1$.
The multinomial chain rule gives
\begin{align}
 A(x,y)&=\sum_{j=0}^{m-1}\ \sum_{|\alpha|=j}
 \frac{D^\alpha A(0)}{\alpha!}(x,y)^\alpha+R_m(x,y),
 \label{eq:componentwise-taylor}\\
 R_m(x,y)&=m\sum_{|\alpha|=m}\frac{(x,y)^\alpha}{\alpha!}
 \int_0^1(1-s)^{m-1}D^\alpha A(sx,sy)\,ds.
 \notag
\end{align}
Continuity bounds the derivatives on a small closed disc, so
$|R_m(x,y)|\le C_m(|x|+|y|)^m=O(|z|^m)$, uniformly in the
angular variable.  If all derivatives of order less than $N$ vanish
at $0$, take $m=N$ to obtain the conclusion.

Otherwise let $k<N$ be the least order of a nonzero derivative.
The nonzero homogeneous Taylor polynomial of degree $k$ has the
unique expression
\[
 P_k(z,\bar z)=\sum_{a+b=k}c_{ab}z^a\bar z^b.
\]
Taking $m=k+1$ in \eqref{eq:componentwise-taylor} yields
\[
 A(re^{i\theta})
 =r^k\sum_{a+b=k}c_{ab}e^{i(a-b)\theta}+O(r^{k+1}),
\]
with a remainder uniform in $\theta$.  Set
\[
 G_k(\theta)=\sum_{a+b=k}c_{ab}e^{i(a-b-N)\theta},
 \qquad F_k(\theta)=G_k(\theta)+\overline{G_k(\theta)}.
\]
It follows that
\begin{equation}\label{eq:componentwise-leading-phase}
 2\operatorname{Re}\!\left(z^{-N}A(z)\right)
 =r^{k-N}\bigl(F_k(\theta)+E(r,\theta)\bigr),
 \qquad \sup_\theta|E(r,\theta)|\le C_Er.
\end{equation}
The \emph{frequencies} of $G_k$ are the pairwise distinct integers
$m_{ab}=a-b-N=2a-k-N\le k-N<0$; those of $\overline{G_k}$
are positive.  Thus $F_k$ cannot vanish identically: Fourier
orthogonality would otherwise give
\[
 0=\frac1{2\pi}\int_0^{2\pi}
 F_k(\theta)e^{-im_{ab}\theta}\,d\theta=c_{ab}
\]
for every $a+b=k$, contradicting $P_k\ne0$.
Moreover, $F_k$ has no zero-frequency term and hence has mean zero.
Being real, continuous and nonzero, it is negative somewhere.
Choose $\theta_0$ and $\delta>0$ with $F_k(\theta_0)=-\delta$.
Along this fixed ray, \eqref{eq:componentwise-leading-phase} becomes
\[
 2\operatorname{Re}\!\left(
 (re^{i\theta_0})^{-N}A(re^{i\theta_0})\right)
 =r^{k-N}\bigl(-\delta+O(r)\bigr)\longrightarrow-\infty.
\]
This contradicts \eqref{eq:componentwise-phase-hypothesis} and proves
$A(z)=O(|z|^N)$.
\end{proof}

Only finite Taylor expansion and Fourier orthogonality enter this
argument; real analyticity is not required.  The lower bound must
hold in every angular direction.  We next establish the geometric
nonemptiness that supplies this condition for divisor masses.

\subsection{Preparation on the punctured disc}
\label{ss:componentwise-punctured}

We first dispose of vertical components and prepare the horizontal
ones away from the exceptional parameter.  For a horizontal component,
we show that its integral-class obstruction vanishes on $\Delta^*$ and
that its restriction there is proper and surjective.
Lemma~\ref{lem:revised-discrete-support} is the key sheaf-theoretic input:
initial vanishing confines the obstruction to a discrete support, and
local freeness on $\Delta^*$ then confines that support to $0$.
This isolates the possible finite-order torsion that must be treated in
Section~\ref{ss:componentwise-mass}.  Surjectivity
ensures the existence of divisors over every nonzero parameter, which
will supply the lower bound in all angular directions required by the
phase lemma.

\begin{lemma}\label{lem:componentwise-punctured}
\begin{enumerate}[$(i)$]
\item Every vertical component of $\mathscr C$ is compact.
\item For every horizontal component $\Lambda$, the restriction
$\Lambda^*:=\Lambda\cap\mu^{-1}(\Delta^*)\to\Delta^*$ is proper and
surjective, and $\Lambda^*$ is dense in $\Lambda$.  Its fixed integral
class $c\in H^2(M,\Z)$ satisfies $o(c)|_{\Delta^*}=0$.
\end{enumerate}
\end{lemma}

\begin{proof}
If $\mu(\Lambda)=\{s\}$, then $\Lambda$ lies in the closed analytic
fiber $\mu^{-1}(s)=\cC_{n-1}(X_s)$.  By its maximality as an
irreducible closed analytic subset of $\mathscr C$, it is also an
irreducible component of $\cC_{n-1}(X_s)$.  Hence $\Lambda$ is compact
by Theorem~\ref{thm:revised-absolute-divisors}.

Suppose now that $\Lambda$ is horizontal.  Lemma~\ref{lem:revised-open-density}
makes its projection open and $\Lambda^*$ dense in $\Lambda$.
Choose a nonempty disc $U\Subset\mu(\Lambda)\cap\Delta^*$.
Under an Ehresmann trivialization, Lemma~\ref{lem:revised-cycle-class} gives
one Poincar\'e-dual class $c$ for every divisor in $\Lambda$.
For each $t\in U$ there is a divisor $Z\in\Lambda_t$, and
$c=c_1(\cO_{X_t}(Z))$ has zero image in $H^2(X_t,\cO_{X_t})$.
The exponential sequence and Lemma~\ref{lem:revised-ddbar-base-change} imply
$o(c)|_U=0$.  Lemma~\ref{lem:revised-discrete-support} now shows that
$S(c)$ is discrete and that every obstruction germ is annihilated by
a finite power of a local coordinate.  Since $R^2\pi_*\cO_{\cX}$ is
locally free, hence torsion-free, on $\Delta^*$, these germs vanish
there.  Thus
\[
 o(c)|_{\Delta^*}=0,\qquad S(c)\subset\{0\}.
\]
Only nonemptiness of $\Lambda_t$ on $U$ is used, not incidence
surjectivity or irreducibility of the divisors themselves.

On a sufficiently small disc $D\Subset\Delta^*$,
Lemma~\ref{lem:revised-local-lift} supplies a total-space line bundle with
Chern class $c$ and a smooth curvature representative $\eta$.
Lemma~\ref{lem:revised-mass-off-S} gives, for every $t\in D$ and
every $Z\in\Lambda_t$,
\[
 \int_ZG_t=\int_{X_t}\eta_t\wedge G_t.
\]
The right-hand side is smooth
in $t$ and independent of $Z$.  A finite cover of any compact
$K\Subset\Delta^*$ by smaller discs with compact closures in such
$D$ gives a uniform mass bound on $\Lambda\cap\mu^{-1}(K)$.
Comparison with a fixed Hermitian metric gives a uniform volume bound.
The cycles have supports in $\pi^{-1}(K)$, and their parameter set is
closed, so Lemma~\ref{lem:revised-bishop-cycle} makes this set compact.
Thus $\Lambda^*\to\Delta^*$ is proper.  Its image is analytic by
Remmert's Theorem~\ref{thm:revised-remmert} and contains $U$;
hence it is all of $\Delta^*$.
\end{proof}

\subsection{The mass estimate on an arbitrary component}
\label{ss:componentwise-mass}

The remaining issue is to bound divisor masses uniformly near the
exceptional parameter.  Lemma~\ref{lem:revised-discrete-support} provides
the crucial link between the sheaf-theoretic obstruction and the analytic
estimate: it supplies a finite exponent $N$ with $t^No(c)=0$ near $0$,
without assuming torsion-freeness of $R^2\pi_*\cO_{\cX}$.
The case $N=0$ recovers the unobstructed situation underlying
Rao--Tsai's torsion-freeness method; the phase estimate handles $N>0$.
We first turn this annihilation relation into a smooth primitive on the
total space, and then express
the mass as a smooth term plus a term to which the phase lemma applies.
Nonnegativity of the masses yields the required lower bound and hence
the uniform estimate.  The argument uses only the fixed integral class
of the component and does not require a good filling.

\begin{lemma}[A primitive from finite-order torsion]
\label{lem:componentwise-torsion-primitive}
Let $c\in H^2(M,\Z)$, let $D_0\subset\Delta$ be a disc about $0$,
and suppose $t^No(c)=0$ on $D_0$ for an integer $N\ge0$.
There are a real smooth $d$-closed two-form $\alpha$ representing
the real image of the corresponding total-space class and a smooth
absolute $(0,1)$-form $u$ on $\cX_{D_0}$ such that
\begin{equation}\label{eq:componentwise-absolute-primitive}
 \bar\partial_{\cX}u=t^N\alpha^{0,2}_{\cX}.
\end{equation}
In particular, $\bar\partial_tu_t=t^N\alpha_t^{0,2}$ on each fiber.
\end{lemma}

\begin{proof}
An Ehresmann trivialization identifies $c$ with a total-space
topological class on $\cX_{D_0}$.  Choose a real smooth closed
representative $\alpha$ of its real image.  Coherence, Cartan's
Theorem~\ref{thm:revised-cartan}, and the Leray spectral sequence give
\[
 H^2(\cX_{D_0},\cO_{\cX_{D_0}})
 \simeq H^0(D_0,R^2\pi_*\cO_{\cX}).
\]
Under the Dolbeault resolution, $o(c)$ corresponds to the class of
$\alpha^{0,2}_{\cX}$, which is $\bar\partial_{\cX}$-closed since
$d\alpha=0$.  Multiplication by $t^N$ on the right corresponds to
multiplication by the holomorphic function $\pi^*t^N$ on the left.
Thus $t^No(c)=0$ says that $t^N\alpha^{0,2}_{\cX}$ has zero
Dolbeault class.  The fine Dolbeault resolution supplies a global
smooth $u$ satisfying \eqref{eq:componentwise-absolute-primitive},
also on the noncompact manifold $\cX_{D_0}$.  Restriction to a fiber
commutes with $\bar\partial$, giving the final assertion.
\end{proof}

In the following estimate the integral is taken over a divisor $Z$
parametrized by $\Lambda$, not over the parameter space $\Lambda$.
The estimate includes reducible and nonreduced divisors.

\begin{prop}[Componentwise mass bound]\label{prop:componentwise-mass-bound}
Let $\Lambda$ be any irreducible component of $\mathscr C$.
For every compact set $K\Subset\Delta$, there is a constant
$C_{\Lambda,K}$ such that
\begin{equation}\label{eq:componentwise-mass-bound}
 0\le\int_ZG_t\le C_{\Lambda,K}
 \qquad(t\in K,\ Z\in\Lambda_t).
\end{equation}
No good-filling assumption on $\Lambda$ is required.
\end{prop}

\begin{proof}
If $\Lambda$ is vertical, it is compact by
Lemma~\ref{lem:componentwise-punctured}.  Continuity of cycle
integration, Lemma~\ref{lem:revised-absolute-G}, gives the bound.
Hence assume $\Lambda$ is horizontal, and let $c$ be its fixed
integral class.  Lemma~\ref{lem:componentwise-punctured} gives the bound on
compact subsets of $\Delta^*$ and nonemptiness of $\Lambda_t$ for
every $t\ne0$.  It remains to control the masses near $0$.

Since $S(c)\subset\{0\}$, Lemma~\ref{lem:revised-discrete-support}
gives a disc $D_0$ about $0$ and the finite exponent $N=N(0)\ge0$
with $t^No(c)=0$ on $D_0$.  Here $N=0$ precisely when the obstruction
vanishes near $0$.  This finite-order control, rather than vanishing
of the obstruction itself, is the essential input for the following
construction.
Lemma~\ref{lem:componentwise-torsion-primitive} supplies a real
smooth $d$-closed two-form $\alpha$
representing the real image of $c$ and a smooth absolute $(0,1)$-form
$u$ on $\cX_{D_0}$ such that
\[
 \bar\partial_{\cX}u=t^N\alpha^{0,2}_{\cX}.
\]
For $t\ne0$, put
\[
 \beta_t=t^{-N}u_t+\bar t^{-N}\bar u_t,
 \qquad \eta_t=\alpha_t-d_t\beta_t.
\]
The form $\eta_t$ is real and $d_t$-closed.  The fiberwise primitive
identity and its conjugate give
\[
 (d_t\beta_t)^{0,2}=t^{-N}\bar\partial_tu_t=\alpha_t^{0,2},
 \qquad
 (d_t\beta_t)^{2,0}=\bar t^{-N}\partial_t\bar u_t=\alpha_t^{2,0}.
\]
Thus the two unwanted components cancel, so $\eta_t$ is of type
$(1,1)$.  Subtracting the exact form $d_t\beta_t$ leaves its real
de Rham class unchanged.
For every $Z\in\Lambda_t$, the current $[Z]-\eta_t$ is $d_t$-exact.
Since $X_t$ satisfies the $\partial\bar\partial$-lemma for $t\ne0$,
Lemma~\ref{lem:revised-current-ddbar} and pairing with $G_t$ eliminate this
difference.  Since
$d_t(\beta_t\wedge G_t)=d_t\beta_t\wedge G_t-\beta_t\wedge d_tG_t$,
Stokes' theorem gives
\[
 \int_{X_t}\eta_t\wedge G_t
 =\int_{X_t}\alpha_t\wedge G_t-\int_{X_t}\beta_t\wedge d_tG_t.
\]
Only $u_t\wedge\partial_tG_t$ and its conjugate have top bidegree
in the last pairing.  Consequently
\begin{align}
 \int_ZG_t
 &=\int_{X_t}\eta_t\wedge G_t
   =B(t)+2\operatorname{Re}\!\left(t^{-N}A(t)\right),
   \label{eq:componentwise-phase-volume}\\
 B(t)&:=\int_{X_t}\alpha_t\wedge G_t,
 \qquad A(t):=-\int_{X_t}u_t\wedge\partial_tG_t.
 \notag
\end{align}
Here reality of $G_t$ gives
$\int_{X_t}\bar u_t\wedge\bar\partial_tG_t
=\overline{\int_{X_t}u_t\wedge\partial_tG_t}$, explaining the
factor $2$ and the real part.  Smooth fiber integration shows that
both $A$ and $B$ are smooth across $0$,
$B$ is real, and the expression is independent of $Z$ within the
fixed component and fiber.

Choose a disc $0\in D\Subset D_0$.  For every $t\in D\setminus\{0\}$,
Lemma~\ref{lem:componentwise-punctured} supplies an effective divisor
in $\Lambda_t$.  Thus \eqref{eq:componentwise-phase-volume} is
nonnegative in every angular direction.  Boundedness of $B$ gives
\[
 \operatorname{Re}\!\left(t^{-N}A(t)\right)
 \ge -\tfrac12\sup_{\overline D}|B|.
\]
For $N\ge1$, Lemma~\ref{lem:componentwise-phase} implies $A(t)=O(|t|^N)$,
so \eqref{eq:componentwise-phase-volume} is bounded above near $0$.
For $N=0$, boundedness follows directly from smoothness of $A$ and
$B$.  This proves a uniform bound for all divisors of $\Lambda$
over a sufficiently small punctured disc.  No estimate on the full
form $t^{-N}u_t$ is needed.

Density of $\Lambda^*$ and continuity of cycle integration extend
the same bound to every member of $\Lambda_0$.  Choose $r>0$ small
enough that the bound holds for $|t|\le r$.  For any compact
$K\Subset\Delta$, the remaining parameters
$K\cap\{|t|\ge r/2\}$ form a compact subset of $\Delta^*$,
where Lemma~\ref{lem:componentwise-punctured} provides a mass bound.  Combining the
two bounds proves \eqref{eq:componentwise-mass-bound} on the original
base $\Delta$.
\end{proof}

\subsection{Properness and the main theorem}
\label{ss:componentwise-conclusion}

We now convert the mass estimate into properness of every component
over the entire disc by Bishop compactness.  Barlet's local
algebraic-dimension theorem then gives Theorem~\ref{thm:main}.
We also explain the alternative good-filling conclusion and finish
with the case of uncountably many Moishezon fibers.

\begin{prop}[Properness of every component]
\label{prop:componentwise-proper-main}
Every irreducible component of $\mathscr C$ is proper over $\Delta$.
Every horizontal component maps surjectively onto $\Delta$.
\end{prop}

\begin{proof}
Fix an irreducible component $\Lambda$ of $\mathscr C$ and a compact set
$K\Subset\Delta$.  Proposition~\ref{prop:componentwise-mass-bound}
bounds the $G_t$-masses of all its divisors over $K$.  Since $G_t$ is
smooth and strictly positive, comparison on the compact set
$\pi^{-1}(K)$ with a fixed Hermitian metric bounds their volumes.
The set $\Lambda\cap\mu^{-1}(K)$ is closed, so
Lemma~\ref{lem:revised-bishop-cycle} makes it compact.  This is properness
over $\Delta$.  For horizontal $\Lambda$, its image is closed and
contains $\Delta^*$ by Lemma~\ref{lem:componentwise-punctured};
it therefore equals $\Delta$.
\end{proof}

\begin{proof}[Proof of Theorem~\ref{thm:main}]
The proper holomorphic submersion $\pi:\cX\to\Delta$ is surjective
and $n$-equidimensional, with irreducible total space and base.
Proposition~\ref{prop:componentwise-proper-main} verifies properness
of every irreducible component of its relative divisor space.
Theorem~\ref{thm:revised-barlet-local}, namely
\cite[Proposition~4.0.8]{Ba15}, now gives an integer
$a_{\mathrm{gen}}\in\{0,\ldots,n\}$ and a countable union $\Sigma$
of proper closed analytic subsets of $\Delta$ such that
\[
 a(X_t)=a_{\mathrm{gen}}\quad(t\in\Delta\setminus\Sigma),
 \qquad a(X_t)\ge a_{\mathrm{gen}}\quad(t\in\Delta).
\]
As $\Sigma$ is countable, every nonempty punctured neighbourhood of
every $s\in\Delta$ meets $\Delta\setminus\Sigma$.  Consequently
\[
 \inf_{0<|t-s|<\varepsilon}a(X_t)=a_{\mathrm{gen}}\le a(X_s)
\]
whenever the indicated disc is contained in $\Delta$.
This also proves \eqref{eq:main-inequality}.  No separate argument
is needed when $a_{\mathrm{gen}}=0$.
\end{proof}

\begin{rem}[The good-filling route]
\label{rem:componentwise-good-filling-route}
One can alternatively follow a selected local good filling rather
than use Barlet's componentwise criterion at the last step.
The local algebraic-dimension
theorem on $\Delta^*$ first gives its common very general value
$a_{\mathrm{gen}}$.  If this value is positive,
Proposition~\ref{prop:revised-local-filling} supplies $\Gamma_U$ whose
generic Kodaira-image dimension is $a_U\ge a_{\mathrm{gen}}$.
Lemma~\ref{lem:revised-global-component} embeds it as a local component of
a full global component $\Gamma$.
Proposition~\ref{prop:componentwise-proper-main} makes $\Gamma$
proper over $\Delta$.

To recover a good filling, consider the addition map for pairs of
nonzero effective relative cycles whose sum belongs to $\Gamma$.
Locally over $\Gamma$ this map is proper: the masses of both summands
are bounded by that of the sum, Bishop compactness gives convergent
subsequences, and the monotonicity inequality gives a positive lower
bound for nonzero cycle masses over a compact part of the base.
Neither summand can therefore disappear in the limit.  Remmert's
Theorem~\ref{thm:revised-remmert} makes the reducible-cycle locus
closed analytic.
Since $\Gamma_U$ supplies irreducible members on a nonempty
open subset, the generic member of $\Gamma$ is irreducible.
The incidence image is closed analytic by properness and Remmert's
Theorem~\ref{thm:revised-remmert}; it contains $\cX_U$, hence equals $\cX$.
Thus $\Gamma$ is a good filling.

One must also retain the local Kodaira-image dimension.  On the
general locus of $U$, write $d=\dim\Gamma_{U,t}$ and choose the
common open subset $\Omega_t$ where the local and global parameter
spaces agree.  Its complement $R_t$ in $\Gamma_{U,t}$ has dimension
at most $d-1$.  The incidence graph over $R_t$ has dimension at most
$d+n-2$, so over a general $x\in X_t$ it contains no component of
the pure $(d-1)$-dimensional local incidence cycle.  That cycle is
therefore recovered, with multiplicities, from its restriction to
$\Omega_t$.  Barlet's relative Kodaira construction
\cite[\S~2.2]{Ba15} identifies the local and global incidence cycles
on this open set.  Equality of global incidence cycles implies
equality of local ones, giving generic image dimension $q\ge a_U$.
Proposition~\ref{prop:revised-barlet-endgame} concludes
$a(X_0)\ge q\ge a_U\ge a_{\mathrm{gen}}$; the case
$a_{\mathrm{gen}}=0$ is immediate.  The proof above retains the
same phase calculation but uses componentwise properness to invoke
Barlet's theorem directly, so neither these good-filling checks nor
the local--global Kodaira comparison need be repeated.
\end{rem}

\begin{proof}[Proof of Theorem~\ref{Moishezon}]
Let $h:=\min_{t\in\Delta}h^1(X_t,\cO_{X_t})$.  By Grauert's
semicontinuity Theorem~\ref{thm:revised-grauert-tools}, the set
\[
 E:=\{t\in\Delta:h^1(X_t,\cO_{X_t})\ge h+1\}
\]
is a proper closed analytic subset of $\Delta$, hence discrete.
At every $s\in\Delta\setminus E$, the Hodge number $h^{0,1}(X_t)$ is
locally constant.  Apply \cite[Theorem~1.4]{RT21} after an automorphism
of the unit disc sending $s$ to $0$.  The theorem requires uncountably
many Moishezon parameters in the whole disc, not necessarily accumulating
at the reference point, so it gives that $X_s$ is Moishezon.

For each $p\in E$, choose a disc $D\Subset\Delta$ with $D\cap E=\{p\}$.
Every fiber over $D\setminus\{p\}$ is Moishezon.  Applying
Corollary~\ref{cor:moishezon} to this restricted family, after rescaling
and recentering the base, shows that $X_p$ is Moishezon as well.
\end{proof}

\section{Scope, necessity, and test families}
\label{s:scope-examples}

The final section tests the conclusion against explicit families.
The Fujiki--Pontecorvo examples explain the role of the cohomological
hypothesis and the difference between a Chern class and a chosen line
bundle.  The K3 families then distinguish the infimum inequality proved
here from ordinary upper semicontinuity, even when every fiber is
K\"ahler.

We distinguish throughout this section between \emph{ordinary upper
semicontinuity} at a parameter $t_0$, expressed by
\begin{equation}\label{eq:limsup-usc}
 a(X_{t_0})\geq\limsup_{t\to t_0}a(X_t),
\end{equation}
and the weaker \emph{punctured-neighbourhood infimum inequality}
\begin{equation}\label{eq:infimum-usc}
 a(X_{t_0})\geq
 \inf_{0<|t-t_0|<\varepsilon}a(X_t).
\end{equation}
Inequality~\eqref{eq:limsup-usc} may fail in smooth families all of whose fibers are
K\"ahler.  Theorem~\ref{thm:main} asserts that the very general value is
the global minimum and, consequently, gives the second inequality.
It asserts neither the first inequality nor upper semicontinuity for
the ordinary analytic Zariski topology.

\subsection{The Fujiki--Pontecorvo jump}
\label{ss:FP-jump}

The Fujiki--Pontecorvo family shows that smoothness alone does not prevent
a downward jump in algebraic dimension.

Let $S$ be a half Inoue surface with $b_2(S)=m$ and let $C$ be its unique
twisted anticanonical curve.  Fujiki--Pontecorvo study the Kuranishi
family of the pair
\[
 g:(\cS,\cC)\longrightarrow T.
\]
The base $T$ is smooth of dimension $m$ and contains a normal-crossings
divisor $A=\bigcup_{j=1}^m A_j$.  The central fiber is $S_0=S$.
For $t\in A\setminus\{0\}$, the fiber is a
blown-up half Inoue surface; for $t\notin A$, it is a blown-up elliptic
diagonal Hopf surface.  Moreover,
\begin{equation*}
 a(S_t)=0\quad(t\in A),\qquad a(S_t)=1\quad(t\notin A).
\end{equation*}
These statements follow from
\cite[Theorem~1.1, Corollary~1.2, and Proposition~1.3]{FP10}.

For $m=1$, after shrinking, $T=\Delta$ and $A=\{0\}$, so
\cite[Example following Proposition~1.3]{FP10} gives
\begin{equation}\label{eq:FP-disk}
 a(S_0)=0<1=\inf_{t\in\Delta^*}a(S_t).
\end{equation}
Thus both \eqref{eq:limsup-usc} and \eqref{eq:infimum-usc} fail at the
central parameter.  However, every fiber is a class VII surface with
$b_1=1$, whereas the $\partial\bar\partial$-lemma implies
\[
 H^1_{\mathrm{dR}}(X,\C)=H^{1,0}_{\bar\partial}(X)
 \oplus H^{0,1}_{\bar\partial}(X),
 \qquad h^{1,0}=h^{0,1},
\]
and hence even $b_1=2h^{0,1}$.  Thus no fiber satisfies that lemma, and
\eqref{eq:FP-disk} shows why the cohomological hypothesis of
Theorem~\ref{thm:main} cannot be dropped.  The punctured fibers are not
Moishezon, since $a(S_t)=1<2$.

Fujiki--Pontecorvo also record higher-dimensional variants in
\cite[Remark~1.1]{FP10}.  The products $S_t^n$ and the Douady spaces
$S_t^{[n]}$ have dimension $2n$ and algebraic dimension $0$ on $A$ and
$n$ off $A$; none is Moishezon.

For the products, the $\partial\bar\partial$-lemma would descend along the
holomorphic retract $S_t^n\rightleftarrows S_t$, so it still fails.  For
the Douady spaces,
$b_1(S_t^{[n]})=b_1(S_t)=1$ by the G\"ottsche formula in
\cite[Theorem~5.2.1.(2)]{DM00}, so it also fails.  The downward jump can
therefore be arbitrarily large outside our hypotheses.  Every punctured
fiber has algebraic dimension $n$ only when $m=1$, or on a disc meeting
$A$ only at its central point.

\subsection{The two Fujiki--Pontecorvo line bundles}
\label{ss:FP-bundles}

This family also distinguishes extension of a specified punctured line
bundle from realization of its Chern class on the total space.

In the one-parameter case, $S_t$ for $t\ne0$ has a canonical elliptic
fibration $f_t:S_t\to\mathbb P^1$.  A fiber decomposes as
\begin{equation*}
 F_t=C_t+E_t,
\end{equation*}
where $E_t$ is the unique $(-1)$-curve.  Consider
\begin{equation}\label{eq:FP-two-bundles}
 \cF_t:=\cO_{S_t}(F_t)=f_t^*\cO_{\mathbb P^1}(1),
 \qquad \cE_t:=\cO_{S_t}(E_t).
\end{equation}
Fujiki--Pontecorvo prove in \cite[Remark~3.2.(2)]{FP10} that neither
prescribed punctured relative bundle extends across the origin.
Indeed, since $h^0(S_t,\cF_t)=2$, an extension $\cF$ would satisfy
$h^0(S_0,\cF_0)\ge2$ by upper semicontinuity.  The quotient of two
independent sections would be a nonconstant meromorphic function on
$S_0$, contrary to $a(S_0)=0$.

The curve $C_t$ extends in the Kuranishi family, and
\[
 \cF_t\simeq\cO_{S_t}(C_t)\otimes\cE_t;
\]
so an extension of $\cE_t$ would also extend $\cF_t$, a contradiction.

This concerns only the representatives in \eqref{eq:FP-two-bundles}, not
every line bundle with the same fiberwise $c_1$; it therefore does not
establish $o(c)\ne0$.  Under its Hodge hypothesis, our theorem would yield
a good filling and $a(X_0)\ge1$, not extension of either prescribed bundle.

Neither bundle in \eqref{eq:FP-two-bundles} is big: $\cF_t$ has Iitaka
dimension one because it defines the elliptic pencil, and a big line bundle
on the surface $S_t$ would make $S_t$ Moishezon, contrary to
$a(S_t)=1<2$.

\subsection{Huybrechts' twistor family}
\label{ss:twistor-family}

Twistor families distinguish the two semicontinuity inequalities even
when all fibers are K\"ahler.

Let $(S,L)$ be a polarized K3 surface with $\rho(S)=1$, and choose the
Ricci-flat K\"ahler metric whose K\"ahler class is proportional to
$c_1(L)$.  Consider the associated twistor family
$\mathscr S\to\mathbb P^1$.

Every fiber, including those on the real equator, is a K3 surface and
hence K\"ahler by Y.-T. Siu's theorem \cite[Main Theorem]{Si83}; all satisfy the
$\partial\bar\partial$-lemma.

For a fixed rational class $\gamma$, type $(1,1)$ means orthogonality to
the moving period line.  Unless this holds identically, it cuts the
twistor conic in at most two points.

Write $T(S)=\operatorname{NS}(S)^\perp\subset H^2(S,\Z)$ for the
\emph{transcendental lattice}.
Classes orthogonal to the entire positive twistor three-plane remain of
type $(1,1)$ on every fiber.  Under the rank-one hypothesis no nonzero
rational class does so: it would lie in $T(S)_\Q$ and be of type
$(1,1)$ on the original fiber, contradicting $\rho(S)=1$.

The rational Hodge structure $T(S)_\Q$ is irreducible and has signature
$(2,19)$.  D. Huybrechts' description therefore applies: the locus where an
additional rational $(1,1)$-class occurs is countable and dense, by
\cite[Introduction and Proposition~3.2]{Hu22}; see also
\cite[\S~1.3]{Hu23}.

A very general twistor fiber has Picard number zero and hence algebraic
dimension zero, whereas the projective fibers are countable and dense.
Thus, on a small coordinate disc about a very general parameter $t_0$,
\begin{equation*}
 a(\mathscr S_{t_0})=0
 <2=\limsup_{t\to t_0}a(\mathscr S_t),
 \qquad
 \inf_{0<|t-t_0|<\varepsilon}a(\mathscr S_t)=0.
\end{equation*}
Thus ordinary upper semicontinuity \eqref{eq:limsup-usc} fails although
every fiber is K\"ahler, while the weaker inequality
\eqref{eq:infimum-usc} holds with equality.  The dense projective fibers
do not give a stronger infimum bound: their polarizations arise from
different integral classes and do not form one relative good filling.

\subsection{Huybrechts' Brauer family}
\label{ss:brauer-family}

In the Brauer family, a persistent genus-one fibration supplies an
explicit proper good filling and a common line bundle, even on
nonalgebraic fibers.

Let $p:S\to\mathbb P^1$ be a projective elliptic K3 surface with a section.
Huybrechts reviews the Brauer family
\begin{equation*}
 \varpi:\mathscr S\longrightarrow H^{0,2}(S)\simeq\C
\end{equation*}
together with a relative genus-one fibration
\begin{equation*}
 g:\mathscr S\longrightarrow\mathbb P^1\times\C;
\end{equation*}
see \cite[\S~1.2, equation~(1.1), and Remark~1.7]{Hu23}.  A detailed
construction through the analytic Tate--Shafarevich group is given in
\cite[Chapter~1.5]{FM94}.

Let $\alpha_t\in\operatorname{Br}^{\mathrm{an}}(S)$ be the image of
$t\in H^{0,2}(S)$ under the natural map to the analytic Brauer group.
By \cite[Corollary~3.3]{Hu23},
\[
 S_t\ \text{is algebraic}
 \quad\Longleftrightarrow\quad
 t\in T(S)\otimes_{\Z}\Q
 \quad\Longleftrightarrow\quad
 \alpha_t\in\operatorname{Br}(S)
 \subset\operatorname{Br}^{\mathrm{an}}(S),
\]
where $T(S)\otimes\Q$ is viewed in $H^{0,2}(S)$ by the natural projection.

Every fiber is a K3 surface, hence K\"ahler by \cite[Main Theorem]{Si83}, and has
$a(S_t)\ge1$ by the genus-one fibration $g_t$.  Since a compact K\"ahler
Moishezon manifold is projective by \cite{Mo67},
\begin{equation*}
 a(S_t)=
 \begin{cases}
 2,&t\in T(S)\otimes_{\Z}\Q,\\
 1,&t\notin T(S)\otimes_{\Z}\Q.
 \end{cases}
\end{equation*}

To exhibit the good filling, for $(b,t)\in\mathbb P^1\times\Delta$ let
$F_{b,t}=g^{-1}(b,t)\subset S_t$.  The classifying map
\[
 \Phi:\mathbb P^1\times\Delta\longrightarrow
 \cC_1(\mathscr S_\Delta/\Delta),
 \qquad(b,t)\longmapsto[F_{b,t}],
\]
is holomorphic by Theorem~\ref{thm:revised-cycle-universal}.$(i)$--$(ii)$,
since $g$ is proper and flat and its fundamental fiber cycles form an
analytic family.

Properness of $\mathbb P^1\times\Delta\to\Delta$ makes $\Phi$ proper over
$\Delta$.  By Theorem~\ref{thm:revised-remmert}, its image $\Gamma$ is
irreducible, analytic, and $\Delta$-proper.  Its generic member is
irreducible and its incidence graph covers $\mathscr S_\Delta$, so
$\Gamma$ is a good filling.

Its relative Kodaira map $g$ has one-dimensional image on every fiber,
and the common fiber class is realized by
\begin{equation*}
 \mathscr L=g^*\operatorname{pr}_1^*\cO_{\mathbb P^1}(1),
\end{equation*}
so its obstruction support is empty.  Barlet's endgame gives
\[
 a(S_0)\ge1=\inf_{t\in\Delta^*}a(S_t).
\]

The image of $T(S)\otimes\Q$ in $H^{0,2}(S)\simeq\C$ is countable and
dense: rational points are dense in $T(S)\otimes\R$, whose projection
onto $H^{0,2}(S)$ is real-linearly surjective.  Its complement is also
dense, so at every nonalgebraic $t_0$ and for every small $\varepsilon>0$,
\begin{equation*}
 a(S_{t_0})=1
 <2=\limsup_{t\to t_0}a(S_t),
 \qquad
 \inf_{0<|t-t_0|<\varepsilon}a(S_t)=1.
\end{equation*}
Thus \eqref{eq:limsup-usc} fails while \eqref{eq:infimum-usc} holds with
equality.  At $0$, the original projective fiber $S$ has algebraic
dimension two, equal to the limsup, whereas the punctured infimum is one.
The dense algebraic parameters do not
produce a relative big class: their polarizations need not arise from one
fixed integral cohomology class.

\appendix
\section[Applications to Moishezon families]
{Applications to Moishezon families\space
{(Joint with Mu-Lin Li, Kai Wang, and Mengjiao Wang)}}
\label{s:applications}

We give several applications of Theorem~\ref{Moishezon} to the deformation invariance of plurigenera, the structure of families, projective loci, and cycle-space theory.

By Theorem~\ref{Moishezon} and the bimeromorphic embedding proved in \cite{RT21,RT22}, a smooth family with uncountably many Moishezon fibers is \emph{Moishezon}, i.e., bimeromorphic to a projective family over the unit disc.
\subsection{Deformation invariance of plurigenera}
The first application to a generically fiberwise Moishezon family is the deformation invariance of plurigenera.
For any positive integer $m$,
the \emph{$m$-genus} $P_{m}(X)$ of a compact complex manifold $X$ is defined by
$$P_{m}(X):= \dim_{\mathbb{C}}H^{0}(X,K_{X}^{\otimes{m}}),$$
where $K_{X}$ is the canonical line bundle of $X$.
It is a bimeromorphic invariant, and the \emph{$m$-genus of an arbitrary compact complex variety} is thus defined as the $m$-genus of any of its nonsingular models.
Similarly, a compact complex variety $M$ is \emph{of general type} if the Kodaira--Iitaka dimension of the canonical line bundle of any of its nonsingular models equals the complex dimension of $M$.
\begin{prop}[{\cite[Main Theorem 1.2]{RT22}}]\label{genus}
Let $\pi: \mathcal{X}\rightarrow\Delta$ be a family such that at least one of the following conditions holds:
\begin{enumerate}[$(i)$]
\item the family $\pi$ is smooth and every fiber $X_{t}$ for $t\in\Delta$ is Moishezon;
\item each (not necessarily projective) compact variety $X_{t}$ for $t\in\Delta$ has only canonical singularities, and uncountably many fibers are of general type.
\end{enumerate}
Then the \emph{$m$-genus} $P_{m}(X_{t})$ for each positive integer $m$ is independent of $t\in\Delta$.
\end{prop}

\begin{cor}
Let $\pi: \cX\rightarrow\Delta$ be a smooth family with uncountably many Moishezon fibers.  Then the $m$-genus $P_{m}(X_{t})$ is independent of $t\in\Delta$ for any positive integer $m$.
\end{cor}
\begin{proof}
This is a direct corollary of Theorem \ref{Moishezon} and Proposition \ref{genus}.
\end{proof}
\begin{rem}
A (possibly incomplete) overview of progress on the deformation invariance of plurigenera can be found in \cite[Introduction]{RT22}.
\end{rem}

\subsection{Pseudo-projective families and projective loci}
Following \cite[Definition 4.28]{RT22}, a family $\pi: \mathcal{X}\rightarrow\Delta$ is called \emph{pseudo-projective} if there exists a bimeromorphic morphism
$$\tau: \mathcal{Y}\rightarrow \cX$$
over $\Delta$, where $\rho: \mathcal{Y}\rightarrow \Delta$ is a projective family with $\mathcal{Y}$ a complex manifold, such that $\tau$ induces a biholomorphism from $\rho^{-1}(U)$ to $\pi^{-1}(U)$ for some dense Zariski open subset $U$ of $\Delta$.

As the second application, one has:
\begin{cor}\label{pp}
Let $\pi: \cX\rightarrow\Delta$ be a smooth family with uncountably many projective fibers.  Then  $\pi$ is pseudo-projective.
\end{cor}
\begin{proof} The proof is implicitly contained in \cite[Lemma 4.25, Remark 4.26]{RT22}. In fact, all fibers are Moishezon by Theorem \ref{Moishezon} and according to the argument in Corollary \ref{main theorem1} below or \cite{RT21,RT22}, it is not difficult to find a (global) line bundle $\mathcal{L}$ on $\mathcal{X}$ such that its restriction to any general fiber $X_{t_0}$ is very ample and $h^{0}(X_t, \mathcal{L}|_{X_t})$ is constant for any $t\in \Delta$ sufficiently close to $t_0$. Then by (the proof of) \cite[Theorem 1.8]{RT22}, the global line bundle $\mathcal{L}$ induces a bimeromorphic map  $\Phi$ from $\mathcal{X}$ to a subvariety  $\mathcal{W}$ of $\mathbb{P}^{m}\times \Delta$ over $\Delta$ for some positive integer $m$, which implies that  $\Phi$ is biholomorphic on $\mathcal{X}\setminus{\mathcal{S}}$ for some proper subvariety $\mathcal{S}$ of $\mathcal{X}$. So $R:=\pi(\mathcal{S})$ is a proper analytic subset of $\Delta$ since the restriction $\Phi_{t_0}:X_{t_0}\rightarrow W_{t_0}$ to the fiber of $\mathcal{W}$ over $t_0$ is biholomorphic.
So $\Phi$ is biholomorphic at least on $\mathcal{X}\setminus{\mathcal{S}}$ by the definition of $\mathcal{S}$, and thus $\mathcal{W}$ is smooth outside ${\pi_{\mathcal{W}}^{-1}(R)}$ for $\pi_{\mathcal{W}}:\mathcal{W}\rightarrow \Delta$ since $\mathcal{X}$ is smooth. Moreover, $\mathcal{W}$ is bimeromorphic to $\mathcal{X}$ and $\pi_{\mathcal{W}}:\mathcal{W}\rightarrow \Delta$ is a projective morphism. By Remmert's elimination of indeterminacy \cite{Re57} and Hironaka's Chow lemma \cite[Corollaries 1 and 2]{Hi75}, one modifies $\mathcal{W}$ to obtain a complex manifold $\mathcal{Y}$ such that the map $\mathcal{Y}\rightarrow {\mathcal{W}}$ over $\Delta$ is a bimeromorphic morphism with $\mathcal{Y}\setminus \pi_{\mathcal{Y}}^{-1}(R)\cong\mathcal{W}\setminus{\pi_{\mathcal{W}}^{-1}(R)}$ as $\mathcal{W}\setminus{\pi_{\mathcal{W}}^{-1}(R)}\cong\mathcal{X}\setminus{\mathcal{S}}$. Hence, $\pi_{\mathcal{Y}}: \mathcal{Y}\rightarrow \Delta$ is the desired projective morphism, and the morphism $\pi$ is pseudo-projective.
\end{proof}

\begin{cor}\label{main theorem1}
Let $\pi: \mathcal{X}\rightarrow\Delta$ be a smooth family. Define the \emph{projective locus} of $\pi$ as
$$\mathcal{P}(\pi):= \{t\in \Delta: \text{$X_{t}$ is projective}\}.$$
Then $\mathcal{P}(\pi)$ is either at most countable or a dense Zariski open subset of $\Delta$.
\end{cor}
It seems not difficult to generalize Corollary \ref{main theorem1} to a proper morphism of complex analytic spaces as in \cite[Section 5]{Ko22b}, where a complete reference list on projective loci is given.
\begin{proof}
This follows directly from Corollary~\ref{pp}. For the reader's convenience, we give a more elementary proof here, also providing further details for the proof of Corollary~\ref{pp}.
Assume that there are uncountably many projective fibers in $\pi$.

\setcounter{step}{0}
\renewcommand{\thestep}{(\Roman{step})} 
\begin{step}\label{step41}
Construct a global line bundle $\cL$ over ${\cX}$ such that $\cL|_{X_{t_0}}$ is ample for some $t_0\in\Delta$.
\end{step}

By the exponential sequence, as in
\eqref{eq:revised-exponential-restriction}, one has the commutative
diagram of exact sequences
\begin{equation*}
\xymatrix@C=0.5cm{
  \cdots \ar[r]^{}
  & H^1(\mathcal{X}, \mathcal{O}^{*}_{\mathcal{X}})\ar[d]_{} \ar[r]^{}
  & H^2(\mathcal{X}, \mathbb{Z}) \ar[d]_{} \ar[r]^{}\ar[d]^{\cong}
  & H^2(\mathcal{X}, \mathcal{O}_{\mathcal{X}})\ar[d]\ar[r]^{} & \cdots \\
   \cdots \ar[r]
  &H^{1}({X_{t}},\mathcal{O}_{X_{t}}^*) \ar[r]^{}
  & H^{2}(X_{t},\mathbb{Z}) \ar[r]
  & H^{2}(X_{t},\mathcal{O}_{{X_{t}}})\ar[r]&\cdots}
\end{equation*}
for any $t\in \Delta$.
Applying the Leray spectral sequence to $\pi$ for the sheaf $\cO_{{\cX}}$, one obtains the exact sequence
$$\xymatrix@C=0.6cm{
   H^1(\Delta,R^1\pi_*\cO_{{\cX}})\ar[r] & H^2({\cX},\cO_{{\cX}}) \ar[r]^{} & H^0(\Delta,R^2\pi_*\cO_{{\cX}}) \ar[r]  & 0. }
$$
Since $\Delta$ is Stein, Theorems~\ref{thm:revised-grauert-tools}.\textup{(i)}
and~\ref{thm:revised-cartan} give
$H^{1}(\Delta,R^1\pi_*\cO_{{\cX}})=0$, and thus
$$H^2({\cX},\cO_{{\cX}})\cong H^0(\Delta,R^2\pi_*\cO_{{\cX}}),$$
as in the Leray argument of
Section~\ref{ss:revised-integral-obstruction}.

Let $c:=c_{1}(A)$ be the first Chern class of an ample line bundle ${A}$ on some fiber $X_{t}$, and denote by $c\in H^{2}(X_{t},\Z)$ the element in $H^2({\cX},\Z)$. The composed morphism
$$H^2({\cX},\Z)\subset H^2({\cX},\C) \to H^0(\Delta,R^2\pi_*\cO_{{\cX}})$$
induces a section $s_{c}\in H^0(\Delta,R^2\pi_*\cO_{{\cX}})$,
namely the obstruction section $o(c)$ of
Section~\ref{ss:revised-integral-obstruction}.
By Theorem~\ref{Moishezon}, every fiber $X_t$ is Moishezon and hence
satisfies the $\partial\bar\partial$-lemma.  The proof of
Lemma~\ref{lem:revised-ddbar-base-change} shows that $h^{0,2}(X_t)$
is constant on $\Delta$.  Thus, by
Theorem~\ref{thm:revised-grauert-tools}.\textup{(iii)},
$R^{2}\pi_{*}{\cO_{\cX}}$ is locally free over $\Delta$ and the map
$$R^2\pi_*\mathcal{O}_{\mathcal{X}}(t):=(R^2\pi_*\mathcal{O}_{\mathcal{X}})_t\otimes \mathbb{C}(t)\rightarrow  H^{2}(X_{t},\mathcal{O}_{{X_{t}}})$$
is bijective for each $t\in\Delta$.
Thus, the zero locus $Z_{c}$ of $s_{c}$ is an analytic subset of $\Delta$. Since there is an uncountable subset $T\subseteq  \Delta$ such that all the fibers $X_{t}$ for $t\in  T$ are projective, one has
\begin{equation*}
\bigcup_{c}Z_{c}\supseteq T,
\end{equation*}
where the union is taken over all integral classes $c_1(L_t)\in H^{2}(X,\Z)$ of ample line bundles $L_t$ on fibers $X_t$ with $t\in T$. Here $X$ is the underlying differentiable manifold of the family's fibers.
Since a countable union of proper analytic subsets is still countable,
there is a class $c$ with $Z_c=\Delta$, so $o(c)=0$.
Lemma~\ref{lem:revised-local-lift}, applied with $W=\Delta$, gives
a global line bundle $\cL$ with $c_1(\cL)=\widetilde c$.
Its restriction to some $X_{t_0}$ has the same first Chern class as an
ample line bundle, and is therefore ample by the Nakai--Moishezon
criterion.

\begin{step}\label{step42}
Show that $\cL|_{X_s}$ is ample outside a proper closed analytic subset
of $\Delta$.
\end{step} 
Following the proof of Grothendieck's theorem
\cite[Th\'eor\`eme~2.1]{Gr61}, we give the details that make the
exceptional set analytic over the whole disc.
Choose $m>0$ such that $\cE:=\cL^{\otimes m}$ restricts to a very ample
line bundle on $X_{t_0}$ and $H^1(X_{t_0},\cE|_{X_{t_0}})=0$.
By Theorem~\ref{thm:revised-grauert-tools}.\textup{(i)},
$\cF:=\pi_*\cE$ is coherent.  Moreover, by
\cite[Corollaire~1.5]{Gr61}, it is locally free near $t_0$ and
\[
 \cF_{t_0}/\mathfrak m_{t_0}\cF_{t_0}
 \simeq H^0(X_{t_0},\cE|_{X_{t_0}}).
\]
By Theorem~\ref{thm:revised-cartan}, applied to
$\mathcal I_{t_0}\cF$ on the Stein disc, global sections of $\cF$
surject onto this fiber.
Thus we can choose $N=h^0(X_{t_0},\cE|_{X_{t_0}})$ global sections
$\sigma_1,\ldots,\sigma_N$ of $\cE$ whose restrictions form a basis
on $X_{t_0}$.  They define a meromorphic map over $\Delta$,
\[
 \phi:\cX\dashrightarrow\mathbb P^{N-1}\times\Delta,
\]
whose restriction $\phi_{t_0}$ is an embedding.  We do not require these
sections to span the complete linear system on any other fiber.

To describe the base points and the failure of fiberwise immersivity
by a single closed analytic set, take local submersion coordinates
$(z_1,\ldots,z_n,t)$ and a holomorphic frame $e$ of $\cE$.
Write $\sigma_i=f_i e$ and set
\[
 J=\begin{pmatrix}
 f_1&\cdots&f_N\\
 \partial f_1/\partial z_1&\cdots&\partial f_N/\partial z_1\\
 \vdots&&\vdots\\
 \partial f_1/\partial z_n&\cdots&\partial f_N/\partial z_n
 \end{pmatrix}.
\]
Changes of frame or submersion coordinates multiply $J$ on the left
by an invertible holomorphic matrix.  Consequently, its
$(n+1)\times(n+1)$ minors define a global closed analytic subset
$Z\subset\cX$.  At a base point its first row vanishes, so the point
belongs to $Z$.  Wherever $f_i\ne0$, elementary row and column
operations give
\[
 \operatorname{rank}J
 =1+\operatorname{rank}
 \bigl(d_{\cX/\Delta}(f_j/f_i)\bigr)_{j\ne i}.
\]
Thus, away from the base locus, $Z$ is exactly the rank-defect locus
of the relative differential of $\phi$.  Since $\phi_{t_0}$ is an
embedding, $Z\cap X_{t_0}=\varnothing$.

By Theorem~\ref{thm:revised-remmert}, $\pi(Z)$ is a proper closed analytic subset of
$\Delta$, hence is discrete.  For $U:=\Delta\setminus\pi(Z)$ and
$s\in U$, the sections have no common zero on $X_s$ and define a
holomorphic immersion $\phi_s:X_s\to\mathbb P^{N-1}$ with
\[
 \cE|_{X_s}\simeq\phi_s^*\cO_{\mathbb P^{N-1}}(1).
\]
The pullback of the Fubini--Study metric has strictly positive curvature,
so $\cE|_{X_s}$, and hence $\cL|_{X_s}$, is ample by the Kodaira
embedding theorem.  In particular, every $X_s$ with $s\in U$ is
projective.  The nonprojective locus is a subset of the closed discrete
set $\pi(Z)$, hence is itself closed analytic in $\Delta$.
Therefore $\mathcal P(\pi)$ is a dense analytic Zariski open subset.
\end{proof}

\subsection{Barlet's theory of relative codimension-one cycle spaces}
As the final (and direct) application of Theorem \ref{Moishezon}, one has:
\begin{cor}
Let $\pi: \cX\rightarrow\Delta$ be a smooth family with uncountably many Moishezon fibers.  Then the restriction of
$$\mu: \mathcal{C}_{n-1}(\cX/\Delta)\rightarrow \Delta$$
to each irreducible component of $\mathcal{C}_{n-1}(\cX/\Delta)$ is proper over $\Delta$.
\end{cor}
\begin{proof}
It follows from Theorem~\ref{Moishezon} and
Lemma~\ref{prop:revised-barlet-proper}.
\end{proof}

\section*{Acknowledgment}
The author is grateful to Professors D. Barlet, F. Campana, Hanlong Fang, Mu-Lin Li, Zhan Li, Yang Shen,  
I-Hsun Tsai, Wei Xia, Song Yang, Xiangdong Yang, and Kang Zuo, as well as
Dr. Jian Chen, Jiawei Feng, Tianzhi Hu, Yi Li, Kai Wang, Mengjiao Wang,
Baidong Zhang and Runze Zhang, for many insightful and helpful discussions on this
topic.  The author also thanks Professors Kefeng Liu, Fangyang Zheng and Xiaojun Huang 
for their continued encouragement and support in his research on
deformation theory and complex geometry.  Part of this work was completed during the author's visit to the School
of Mathematics and Statistics at Hainan University over the Lunar New
Year in 2025.  The author thanks the school for its warm hospitality
and excellent working environment.


The author owes a deep mathematical debt to I-Hsun and Mu-Lin. The author's systematic study of the deformation-limit problem began
during his one-year visit to Professors Chin-Yu Hsiao and Jih-Hsin Cheng since September 2017.  At that stage, the author--I-Hsun were not even aware that a coherent analytic sheaf on a one-dimensional
complex manifold is torsion-free if and only if it is locally free.
This ignorance led us to study torsion sections through their supports, as in
the proof of \cite[Proposition~4.5]{RT21} and also underlies
the treatment of the obstruction $o(c)$ in the present paper. In April 2023, the author resumed work on this topic,  together with Mengjiao Wang, Mu-Lin and
Kai Wang. We had explored numerous
approaches, including the extension of period maps and the
Oka--Grauert principle, etc.
Beginning in late August 2026, much inspired by a recent joint work with Jiawei and Yi Li, the author returned independently to the original
torsion-freeness approach developed with I-Hsun, seeking to use
the differential forms arising from the annihilation relation
$t^No(c)=0$ in the mass estimate of
Proposition~\ref{prop:componentwise-mass-bound}.
During discussions with the author, ChatGPT suggested an initial model
for the phase lemma, developed here as
Lemma~\ref{lem:componentwise-phase}.  These supply the mass bound
sought in Barlet's good-filling approach \cite{Ba15} extends the Rao--Tsai
method beyond the torsion-free case.

\raggedbottom 

\end{document}